\documentclass[11pt,a4paper]{article}
\usepackage{graphicx}
\usepackage{amsmath}
\usepackage{graphicx}
\usepackage{verbatim}
\usepackage{color}
\usepackage{subfigure}
\usepackage{float}
\usepackage{authblk}
\usepackage[colorlinks]{hyperref}
\usepackage{doi}

\usepackage{graphicx}%
\usepackage{multirow}%
\usepackage{amsmath,amssymb,amsfonts}%
\usepackage{amsthm}%
\usepackage{mathrsfs}%
\usepackage[title]{appendix}%
\usepackage{xcolor}%
\usepackage{textcomp}%
\usepackage{manyfoot}%
\usepackage{booktabs}%
\usepackage{algorithm}%
\usepackage{algorithmicx}%
\usepackage{algpseudocode}%
\usepackage{listings}%

\usepackage{textcomp}
\usepackage{eurosym}

\theoremstyle{thmstyleone}%
\newtheorem{theorem}{Theorem}[section]
\newtheorem{proposition}[theorem]{Proposition}%

\theoremstyle{thmstyletwo}%
\newtheorem{remark}[theorem]{Remark}%
\newtheorem{lemma}[theorem]{Lemma}
\newtheorem{corollary}[theorem]{Corollary}
\theoremstyle{thmstylethree}%
\newtheorem{definition}[theorem]{Definition}%
\begin{document}
\title{Counting of lattices containing up to five comparable reducible elements and having nullity up to three}
\author[1]{Balasaheb P. Aware}
\author[2]{Ashok N. Bhavale}

\affil[1]{Department of Mathematics, Modern College of ASC(A), Shivajinagar, Pune, India  (E-mail: aware66pp@gmail.com)}
\affil[2]{Head, Department of Mathematics, Modern College of ASC(A), Shivajinagar, Pune, India  (E-mail: hodmaths@moderncollegepune.edu.in)}

\maketitle
\begin{abstract}In 2020 Bhavale and Waphare introduced the concept of a nullity of a poset as nullity of its cover graph. In 2003 Pawar and Waphare counted all non-isomorphic lattices on n elements and n edges, which are precisely lattices of nullity one. In 2002 Thakare et al. counted all non-isomorphic lattices on n elements containing two reducible elements. In the same paper, Thakare et al. counted lattices on n elements containing up to n+1 edges, which are precisely lattices of nullity up to two. In 2024 Bhavale and Aware counted all non-isomorphic lattices on n elements, containing up to three reducible elements. Recently, Aware and Bhavale counted all non-isomorphic lattices on n elements, containing four comparable reducible elements, and having nullity three. In this paper, we count up to isomorphism the class of all lattices on n elements containing five comparable reducible elements, and having nullity three.

\vskip1em \noindent \textbf{2020 AMS Classification: 06A05; 06A06; 06A07.}


\vskip1em \noindent \textbf{Article type: Research article}
\end{abstract}

\section{Introduction}\label{sec1}

In $1940$ Birkhoff \cite{bib6} raised the open problem, compute for small $n$ all non-isomorphic posets/lattices on a set of $n$ elements. This NP-complete problem was attempted by many authors from all over the world. For more details about counting of lattices see \cite{bib7}, \cite{bib9}, \cite{bib14}.

An element $x$ in a lattice $L$ is $join$-$reducible$($meet$-$reducible$) in $L$ if there exist $y,z \in L$ both distinct from $x$, such that $y\vee z=x (y\wedge z=x)$. An element $x$ in a lattice $L$ is $reducible$ if it is either join-reducible or meet-reducible. $x$ is $join$-$irreducible$($meet$-$irreducible$) if it is not join-reducible(meet-reducible); $x$ is $doubly$ $irreducible$ if it is both join-irreducible and meet-irreducible. The set of all reducible elements in $L$ is denoted by $Red(L)$, and its complement in $L$ is denoted by $Irr(L)$. An element $x$ in a poset $P$ {\it{covers}} an element $y$ in a poset $P$ if $x<y$, and there is no $z$ in $P$ such that $x<z<y$. We denote this fact as $x\prec y$, and we say that $<x,y>$ is {\it{covering}} or an edge. The {\it{height}} of a lattice is the number of coverings in a maximal chain of a lattice. Lattices $L_1$ and $L_2$ are {\it{isomorphic}} (in symbol, $L_1\cong L_2$), and the map $\phi:L_1\to L_2$ is an {\it{isomorphism}} if and only if $\phi$ is one-to-one and onto, and $a \leq b$ in $L_1$ if and only if $\phi(a) \leq \phi(b)$ in $L_2$. Let $L_1$ and $L_2$ be two lattices. By $L_1\oplus L_2(L_1\circ L_2)$ we mean linear sum(vertical sum) of lattices(see \cite{bib19}, \cite{bib20}). If $P$ is a poset then $Irr(P)$ is the set of all elements in $P$ which have at most one lower cover and at most one upper cover in $P$. Bhavale and Waphare \cite{bib5} introduced the concept of a {\it{nullity of a poset}} $P$, denoted by $\eta(P)$, as nullity of its cover graph.  According to Bhavale and Waphare \cite{bib5}, if a dismantlable lattice of nullity $k$ contains $r$ reducible elements then $2\leq r\leq 2k$. Bhavale and Waphare \cite{bib5} also introduced the concept of an {\it{RC-lattice}}, a {\it{basic block}} and a {\it{basic block associated to a lattice}}. 
 
\begin{definition}\cite{bib5}\textnormal{An {\it{RC-lattice}} is a lattice in which all the reducible elements are comparable.}
\end{definition}
\begin{definition}\cite{bib5}\label{basicblock}
\textnormal{A poset $P$ is a {\it{basic block}} if it is one element or $Irr(P) = \emptyset$ or removal of any doubly irreducible element reduces nullity by one}.
\end{definition}

\begin{definition}\cite{bib5}\label{bbas}
\textnormal{$B$ is a $basic$ $block$ $associated$ $to$ $a$ $poset$ $P$ if $B$ is obtained from the basic retract associated to $P$ by successive removal of all the pendant vertices.}
\end{definition}
\begin{theorem}\cite{bib5}\label{redb}
Let $B$ be the basic block associated to a poset $P$. Then $Red(B)=Red(P)$ and $\eta(B)=\eta(P)$.
\end{theorem}
In 1974, Rival \cite{bib18} introduced and studied the class of dismantlable lattices. 
\begin{definition}\cite{bib18}
\textnormal{A finite lattice of order $n$ is called {\it{dismantlable}} if there exist a chain $L_{1} \subset L_{2} \subset \ldots\subset L_{n}(=L)$ of sublattices of $L$ such that $|L_{i}| = i$, for all $i$}.
\end{definition}
Bhavale and Waphare \cite{bib5} proved that every RC-lattice is a dismantlable lattice.
\begin{theorem}\cite{bib14}\label{dac} A finite lattice is dismantlable if and only if it is an adjunct of chains.
\end{theorem}

If a lattice $L$ is an adjunct of chains $C_0,C_1,\ldots, C_k$ then it is written as $L=C_0]_{a_1}^{b_1}C_1]_{a_2}^{b_2}C_2\ldots $\\$C_{k-1}]_{a_k}^{b_k}C_k$, where for $1\leq i\leq k$, $(a_i,b_i)$ is an adjunct pair(see \cite{bib14}).
\begin{theorem}\cite{bib1}\label{bhavale} Let $L$ be a dismantlable lattice, and $C$ be a maximal chain in $L$. Then there exist chains $C_1,C_2,\ldots, C_k$ in $L$ such that $L=C]_{a_1}^{b_1}C_1]_{a_2}^{b_2}C_2\cdots]_{a_k}^{b_k}C_k$.
\end{theorem}
\begin{theorem}\cite{bib5}\label{nulladj}
A dismantlable lattice $L$ containing $n$ elements is of nullity $l$ if and only if $L$ is an adjunct of $l+1$ chains.
\end{theorem}Using Theorem \ref{bhavale} and Theorem \ref{nulladj}, we have the following result.
\begin{corollary}\label{maxchain}
An RC-lattice $L$ containing $n$ elements is of nullity $l$ if and only if $L$ is an adjunct of $l+1$ chains, starting with a maximal chain containing all the reducible elements.
\end{corollary}
In $2003$ Pawar and Waphare \cite{bib11} counted all non-isomorphic lattices with $n$ elements and $n$ edges, which are precisely lattices of nullity one.
\begin{theorem}\label{1.1}\cite{bib11}
$n\geq 4$, $|\mathscr{L}(n;2,1)|=
\begin{cases} 
      \frac{m(m-1)(4m+1)}{6};& if~n=2m+1 \\
      \frac{m(m-1)(4m-5)}{6};&if~n=2m. 
\end{cases}$
\end{theorem}
 Thakare et al. \cite{bib14} generalised this result for arbitrary nullity $k\geq 1$.
\begin{theorem}(\cite{bib14},\cite{bib2})\label{2red}
\textnormal{For} $n\geq 4$ and for $1\leq k\leq n-3$, $|\mathscr{L}(n;2,k)|=\displaystyle\sum_{j=1}^{n-k-2}jP_{n-j-1}^{k+1}$.
\end{theorem}
 Recently, Bhavale and Aware \cite{bib2} counted all non-isomorphic lattices on $n$ elements, containing up to three reducible elements, and having arbitrary nullity $k\geq 2$.

\begin{theorem}\cite{bib2}\label{BA3red} For an integer $n\geq 6$,
$ |\mathscr{L}(n;3,k)|=$\\$\displaystyle\sum_{j=0}^{n-6}\sum_{k=2}^{n-j-4}\sum_{l=1}^{n-j-5}\sum_{i=1}^{n-j-l-4}2(j+1)P^{k}_{n-j-l-i-2}+\sum_{j=0}^{n-6}\sum_{k=3}^{n-j-4}\sum_{r=5}^{n-j-2}\sum_{s=1}^{k-2}\sum_{i=1}^{r-4}2(j+1)P^{s+1}_{r-i-2}P^{k-s}_{n-j-r}+\\ \sum_{j=0}^{n-7}\sum_{k=2}^{n-j-5}\sum_{l=4}^{n-j-3}\sum_{t=1}^{k-1}(j+1)P^{t+1}_{l-2}P^{k-t+1}_{n-j-l-1}+\sum_{j=0}^{n-8}\sum_{k=3}^{n-j-5}\sum_{r=1}^{n-j-7}\sum_{l=4}^{n-j-r-3}\sum_{t=1}^{k-2}(j+1) P^{t+1}_{l-2}P^{k-t}_{n-j-r-l-1}\\+\sum_{j=0}^{n-8}\sum_{k=4}^{n-j-5}\sum_{r=2}^{n-j-7}\sum_{s=2}^{k-2}\sum_{l=4}^{n-j-r-3}
\sum_{t=1}^{k-s-1}(j+1)P^{t+1}_{l-2}P^{k-s-t+1}_{n-j-r-l-1}P^{s}_{r}.$
\end{theorem} 
 
In particular for $k=2,3$,
\begin{corollary}\label{BA3red2} For an integer $n\geq 6$,\\
$\displaystyle |\mathscr{L}(n;3,2)|=\sum_{j=0}^{n-6}\sum_{l=1}^{n-j-5}\sum_{i=1}^{n-j-l-4}2(j+1)P^{2}_{n-j-l-i-2}+\sum_{j=0}^{n-7}\sum_{l=4}^{n-j-3}(j+1)P^{2}_{l-2}P^{2}_{n-j-l-1}$.
\end{corollary}

\begin{corollary}\label{BA3red3} For an integer $n\geq 6$,\\
$|\mathscr{L}(n;3,3)|=\displaystyle\sum_{j=0}^{n-6}\sum_{l=1}^{n-j-5}\sum_{i=1}^{n-j-l-4}2(j+1)P^{3}_{n-j-l-i-2}+\sum_{j=0}^{n-6}\sum_{r=5}^{n-j-2}\sum_{i=1}^{r-4}2(j+1)P^{2}_{r-i-2}P^{2}_{n-j-r}+$\\$\displaystyle \sum_{j=0}^{n-7}\sum_{l=4}^{n-j-3}\sum_{t=1}^{2}(j+1)P^{t+1}_{l-2}P^{4-t}_{n-j-l-1}+\sum_{j=0}^{n-8}\sum_{r=1}^{n-j-7}\sum_{l=4}^{n-j-r-3}(j+1) P^{2}_{l-2}P^{2}_{n-j-r-l-1}$.
\end{corollary}

Thakare et al. \cite{bib14} counted lattices on $n$ elements and having nullity up to two using set theoretic approach. Whereas, Bhavale and Aware \cite{bib3} also counted all non-isomorphic lattices on $n$ elements and having nullity up to two using the concept of {\it{basic block}} and {\it{basic block associated to a lattice}}. Further Bhavale and Aware \cite{bib3} verified that the formulae(lattice having nullity up to two) are equivalent. Moreover, according to our point of view it seems that, the technique used by Bhavale and Aware \cite{bib3} is practically more beneficial for the further counting purpose.
 
\begin{theorem}\label{1.4}\cite{bib3}
For $n\geq 6$, $\displaystyle|\mathscr{L}(n;4,2)| = \sum_{i=0}^{n-6} (i+1) \binom{n-i-2}{4}+ \\
 \sum_{i=0}^{n-7}\sum_{p=1}^{n-i-6}\sum_{l=2}^{n-i-p-4}(i+1)(l-1)P^{2}_{n-i-p-l-2}+ 
 \sum_{i=0}^{n-8}\sum_{m=0}^{n-i-8}\sum_{s=4}^{n-i-m-4} (i+1)(n-i-m-7)P^{2}_{s-2}P^{2}_{n-i-m-s-2}$.
\end{theorem}
 Recently Bhavale and Aware \cite{bib17} also counted all non-isomorphic RC-lattices on n elements, containing four reducible elements, having nullity three.
\begin{theorem}\cite{bib17}
For $n\geq 7$, \\$|\mathscr{L}(n;4,3)|=\scriptsize{\displaystyle\sum_{q=0}^{n-7}(q+1)\bigg(\sum_{s=1}^{n-q-6}\sum_{r=1}^{n-q-s-5} \sum_{l=2}^{n-q-s-r-3}2(n-q-s-r-l-2)P_{l}^2+}$\\$\displaystyle\sum_{p=4}^{n-q-4}\sum_{t=1}^{n-q-p-3}tP_{n-q-p-t-1}^{2}P^2_{p-2}+\scriptsize{\displaystyle\sum_{t=1}^{n-q-7}\sum_{i=2}^{n-q-t-5}(i-1)P^{3}_{n-q-t-i-2}+}$\\$\scriptsize{\displaystyle\sum_{p=1}^{n-q-6}\binom{n-q-p-2}{4}}\scriptsize{+\displaystyle\sum_{t=1}^{n-q-7}\sum_{r=1}^{n-q-t-6}\sum_{l=1}^{n-q-t-r-5}\sum_{i=1}^{n-q-t-r-l-4}7P_{n-q-t-r-l-i-2}^{2}}$\\$\scriptsize{+\displaystyle\sum_{r=0}^{n-q-9}\sum_{p=5}^{n-q-r-4}2P_{p-2}^{3}P_{n-q-p-r-2}^{2}+\displaystyle\sum_{p=4}^{n-q-5}\sum_{l=1}^{n-q-p-4}\sum_{i=1}^{n-q-p-l-3}4P_{p-2}^{2}P_{n-q-p-l-i-1}^2+}$\\$\scriptsize{\displaystyle\sum_{r=1}^{n-q-8}\sum_{q=1}^{n-q-r-7}\sum_{l=4}^{n-q-q-r-3}2P_{l-2}^2P_{n-q-q-r-l-1}^2+}\scriptsize{\displaystyle\sum_{p=7}^{n-q-3}\sum_{l=4}^{p-3}P_{n-q-p-1}^2P_{l-2}^2P_{p-l-1}^2}+$\\$\scriptsize{\displaystyle\sum_{t=1}^{n-q-8}\sum_{m=0}^{n-q-t-8}\sum_{s=4}^{n-q-t-m-4}(n-q-t-m-7)P^{2}_{s-2}P^{2}_{n-q-t-m-s-2}}\bigg)$.
\end{theorem} 
  
Thus counting of all non-isomorphic lattices on $n$ elements containing $r\leq 4$ comparable reducible elements, and having nullity $k\leq 3$ is taken care of. In this paper, we count up to isomorphism the class of all non-isomorphic lattices on $n$ elements, containing exactly five comparable reducible elements, and having nullity three. For this purpose, we need the following two results.
\begin{lemma}\label{oplus}
Let $\mathscr{L}_{1}(p)$ and $\mathscr{L}_{2}(q)$ be some classes of non-isomorphic lattices on $p\geq r\geq 1$ and $q\geq s\geq 1$ elements respectively. If $L_1\in \mathscr{L}_{1}(p)$ and $L_2\in \mathscr{L}_{2}(q)$ are such that $L_1\oplus L_2=L\in \mathscr{L}(n)$, a class of lattices on $n=p+q$ elements then $|\mathscr{L}(n)|=\displaystyle\sum_{p=r}^{n-s}(|\mathscr{L}_{1}(p)|\times|\mathscr{L}_{2}(n-p)|)$.
\end{lemma}
\begin{proof}
As $n=p+q$, for fixed $p$, there are $|\mathscr{L}_{1}(p)|$ non-isomorphic lattices on $p$ elements and $|\mathscr{L}_{2}(n-p)|$ non-isomorphic lattices on $n-p$ elements. Therefore by multiplication principle, for fixed $p$, there are up to isomorphism $|\mathscr{L}_{1}(p)|\times|\mathscr{L}_{2}(n-p)|$ lattices in $\mathscr{L}(n)$. Note that $r\leq p=n-q\leq n-s$, since $q\geq s$. Thus there are $\displaystyle\sum_{p=r}^{n-s}(|\mathscr{L}_{1}(p)|\times|\mathscr{L}_{2}(n-p)|)$ lattices in $\mathscr{L}(n)$ up to isomorphism.
\end{proof}
As vertical sum $L_1\circ L_2$ of lattices $L_1$ and $L_2$ is obtained from the direct sum $L_1\oplus L_2$ of lattices $L_1$ and $L_2$ by identifying the greatest element of $L_1$ with the smallest element of $L_2$. Therefore by Lemma \ref{oplus} we have the following result.
\begin{corollary}\label{circ}
Let $\mathscr{L}_{1}(p)$ and $\mathscr{L}_{2}(q)$ be some classes of non-isomorphic lattices on $p\geq r\geq 1$ and $q\geq s\geq 1$ elements respectively. If $L_1\in \mathscr{L}_{1}(p)$ and $L_2\in \mathscr{L}_{2}(q)$ are such that $L_1\circ L_2=L\in \mathscr{L}(n)$, a class of lattices on $n=p+q-1$ elements then $|\mathscr{L}(n)|=\displaystyle\sum_{p=r}^{n-s+1}(|\mathscr{L}_{1}(p)|\times|\mathscr{L}_{2}(n-p+1)|)$.
\end{corollary}
If $L=C\oplus \textbf{B} \oplus C'$ where $\textbf{B}$ is the {\it{block}}(a lattice in which $0$ and $1$ are reducible elements), and $C$, $C'$ are chains then $\textbf{B}$ is called a {\it{maximal block}}. For the necessary definitions, notation, and terminology see \cite{bib2,bib8,bib15,bib17}.

\section{Counting of RC-lattices containing five reducible elements and having nullity three}\label{sec3}
In this section, we count all non-isomorphic lattices on $n$ elements, containing five comparable reducible elements, and having nullity three. Let $\mathscr{L}(n;r,k)$ be the class of all non-isomorphic RC-lattices on $n$ elements such that every member of it contains $r$ reducible elements, and has nullity $k$. Let $\mathscr{L}(n;r,k,h)$ be the subclass of $\mathscr{L}(n;r,k)$ such that the basic block associated to a member of it is of height $h$. Let $\mathscr{B}(j;r,k)$ be the class of all non-isomorphic maximal blocks on $j$ elements such that every member of it contains $r$ comparable reducible elements and has nullity $k$. Let $\mathscr{B}(j;r,k,h)$ be the subclass of $\mathscr{B}(j;r,k)$ such that the basic block associated to a member of it is of height $h$.
\begin{center}
\unitlength 1mm 
\linethickness{0.4pt}
\ifx\plotpoint\undefined\newsavebox{\plotpoint}\fi 
\begin{picture}(159.448,66.078)(0,0)
\put(48.665,9.276){\circle{1.026}}
\put(48.715,13.288){\circle{1.026}}
\put(48.665,17.3){\circle{1.026}}
\put(48.715,21.272){\circle{1.026}}
\put(48.697,20.867){\line(0,-1){3.124}}
\put(48.697,16.884){\line(0,-1){3.124}}
\put(48.697,12.843){\line(0,-1){3.124}}
\put(48.684,25.211){\circle{1.026}}
\put(48.665,24.806){\line(0,-1){3.124}}
\put(44.923,13.07){\circle{1.026}}
\multiput(48.082,17.363)(-.033642857,-.039183673){98}{\line(0,-1){.039183673}}
\multiput(44.785,12.72)(.034059406,-.033475248){101}{\line(1,0){.034059406}}
\put(48.764,28.952){\circle{1.026}}
\put(48.745,28.547){\line(0,-1){3.124}}
\put(48.733,32.891){\circle{1.026}}
\put(48.714,32.486){\line(0,-1){3.124}}
\put(44.397,25.246){\circle{1.026}}
\multiput(48.298,28.892)(-.040628866,-.033597938){97}{\line(-1,0){.040628866}}
\multiput(44.357,24.898)(.038262136,-.033679612){103}{\line(1,0){.038262136}}
\put(44.282,28.405){\circle{1.026}}
\multiput(44.222,28.114)(.033487603,-.05492562){121}{\line(0,-1){.05492562}}
\multiput(44.323,28.843)(.033478261,.035243478){115}{\line(0,1){.035243478}}
\put(59.949,9.382){\circle{1.026}}
\put(59.999,13.394){\circle{1.026}}
\put(59.949,17.406){\circle{1.026}}
\put(59.999,21.378){\circle{1.026}}
\put(59.98,20.973){\line(0,-1){3.124}}
\put(59.98,16.99){\line(0,-1){3.124}}
\put(59.98,12.949){\line(0,-1){3.124}}
\put(59.968,25.317){\circle{1.026}}
\put(59.949,24.912){\line(0,-1){3.124}}
\put(56.207,13.176){\circle{1.026}}
\multiput(56.069,12.826)(.034059406,-.033475248){101}{\line(1,0){.034059406}}
\put(60.048,29.058){\circle{1.026}}
\put(60.029,28.653){\line(0,-1){3.124}}
\put(60.017,32.997){\circle{1.026}}
\put(59.998,32.592){\line(0,-1){3.124}}
\put(55.565,28.511){\circle{1.026}}
\multiput(55.607,28.949)(.033478261,.035243478){115}{\line(0,1){.035243478}}
\multiput(55.607,28.138)(.04409412,-.03336471){85}{\line(1,0){.04409412}}
\multiput(59.453,21.338)(-.03345,-.07655){100}{\line(0,-1){.07655}}
\put(70.893,9.305){\circle{1.026}}
\put(70.943,13.317){\circle{1.026}}
\put(70.893,17.329){\circle{1.026}}
\put(70.943,21.301){\circle{1.026}}
\put(70.924,20.896){\line(0,-1){3.124}}
\put(70.924,16.913){\line(0,-1){3.124}}
\put(70.924,12.872){\line(0,-1){3.124}}
\put(70.912,25.24){\circle{1.026}}
\put(70.893,24.835){\line(0,-1){3.124}}
\put(67.15,13.099){\circle{1.026}}
\multiput(70.309,17.392)(-.033632653,-.039183673){98}{\line(0,-1){.039183673}}
\multiput(67.013,12.749)(.034059406,-.033475248){101}{\line(1,0){.034059406}}
\put(70.992,28.981){\circle{1.026}}
\put(70.973,28.576){\line(0,-1){3.124}}
\put(70.96,32.92){\circle{1.026}}
\put(70.942,32.515){\line(0,-1){3.124}}
\put(66.602,30.562){\circle{1.026}}
\multiput(66.593,30.109)(.033555556,-.045425926){108}{\line(0,-1){.045425926}}
\multiput(66.593,31.075)(.07286275,.0335098){51}{\line(1,0){.07286275}}
\put(66.602,26.548){\circle{1.026}}
\multiput(70.402,32.784)(-.033707965,-.050646018){113}{\line(0,-1){.050646018}}
\multiput(66.5,26.095)(.0343,-.033730769){130}{\line(1,0){.0343}}
\put(92.458,9.379){\circle{1.026}}
\put(92.508,13.391){\circle{1.026}}
\put(92.458,17.403){\circle{1.026}}
\put(92.508,21.375){\circle{1.026}}
\put(92.489,20.97){\line(0,-1){3.124}}
\put(92.489,16.987){\line(0,-1){3.124}}
\put(92.489,12.946){\line(0,-1){3.124}}
\put(92.477,25.314){\circle{1.026}}
\put(92.458,24.909){\line(0,-1){3.124}}
\put(92.557,29.055){\circle{1.026}}
\put(92.538,28.65){\line(0,-1){3.124}}
\put(92.525,32.994){\circle{1.026}}
\put(92.507,32.589){\line(0,-1){3.124}}
\put(88.074,28.508){\circle{1.026}}
\multiput(88.115,28.946)(.033478261,.035243478){115}{\line(0,1){.035243478}}
\multiput(88.065,28.102)(.04425,-.03361905){84}{\line(1,0){.04425}}
\put(87.982,17.406){\circle{1.026}}
\multiput(91.875,21.561)(-.03612963,-.033722222){108}{\line(-1,0){.03612963}}
\put(87.973,17.919){\line(0,1){0}}
\multiput(87.973,16.953)(.036811321,-.033650943){106}{\line(1,0){.036811321}}
\put(87.888,21.345){\circle{1.026}}
\multiput(91.782,25.426)(-.036820755,-.033660377){106}{\line(-1,0){.036820755}}
\multiput(87.879,20.967)(.033725926,-.083140741){135}{\line(0,-1){.083140741}}
\put(137.684,9.305){\circle{1.026}}
\put(137.734,13.317){\circle{1.026}}
\put(137.684,17.329){\circle{1.026}}
\put(137.734,21.301){\circle{1.026}}
\put(137.715,20.896){\line(0,-1){3.124}}
\put(137.715,16.913){\line(0,-1){3.124}}
\put(137.715,12.872){\line(0,-1){3.124}}
\put(137.703,25.24){\circle{1.026}}
\put(137.684,24.835){\line(0,-1){3.124}}
\put(137.783,28.981){\circle{1.026}}
\put(137.764,28.576){\line(0,-1){3.124}}
\put(137.752,32.92){\circle{1.026}}
\put(137.733,32.515){\line(0,-1){3.124}}
\put(133.199,17.845){\line(0,1){0}}
\put(133.208,19.413){\circle{1.026}}
\multiput(137.659,32.487)(-.033718519,-.093592593){135}{\line(0,-1){.093592593}}
\multiput(133.199,19.035)(.033646552,-.083301724){116}{\line(0,-1){.083301724}}
\put(147.77,13.368){\circle{1.026}}
\put(147.82,17.38){\circle{1.026}}
\put(147.77,21.392){\circle{1.026}}
\put(147.82,25.364){\circle{1.026}}
\put(147.802,24.959){\line(0,-1){3.124}}
\put(147.802,20.976){\line(0,-1){3.124}}
\put(147.802,16.935){\line(0,-1){3.124}}
\put(147.789,29.303){\circle{1.026}}
\put(147.77,28.898){\line(0,-1){3.124}}
\put(147.869,33.044){\circle{1.026}}
\put(147.85,32.639){\line(0,-1){3.124}}
\put(147.838,36.983){\circle{1.026}}
\put(147.819,36.578){\line(0,-1){3.124}}
\put(143.285,21.908){\line(0,1){0}}
\put(143.202,32.46){\circle{1.026}}
\multiput(147.718,36.603)(-.040810811,-.033441441){111}{\line(-1,0){.040810811}}
\multiput(143.188,32.095)(.0485,-.03341463){82}{\line(1,0){.0485}}
\put(147.669,9.374){\circle{1.026}}
\put(147.7,12.941){\line(0,-1){3.124}}
\put(143.202,25.246){\circle{1.026}}
\multiput(147.234,29.208)(-.038433962,-.033716981){106}{\line(-1,0){.038433962}}
\multiput(142.898,24.898)(.042141509,-.033716981){106}{\line(1,0){.042141509}}
\put(142.939,13.473){\circle{1.026}}
\multiput(147.103,17.33)(-.042097087,-.033679612){103}{\line(-1,0){.042097087}}
\multiput(142.767,13.02)(.039669811,-.033707547){106}{\line(1,0){.039669811}}
\put(158.149,13.263){\circle{1.026}}
\put(158.199,17.275){\circle{1.026}}
\put(158.149,21.287){\circle{1.026}}
\put(158.199,25.259){\circle{1.026}}
\put(158.18,24.854){\line(0,-1){3.124}}
\put(158.18,20.871){\line(0,-1){3.124}}
\put(158.18,16.83){\line(0,-1){3.124}}
\put(158.168,29.198){\circle{1.026}}
\put(158.149,28.793){\line(0,-1){3.124}}
\put(158.248,32.939){\circle{1.026}}
\put(158.229,32.534){\line(0,-1){3.124}}
\put(158.217,36.878){\circle{1.026}}
\put(158.198,36.473){\line(0,-1){3.124}}
\put(153.664,21.803){\line(0,1){0}}
\put(153.58,32.355){\circle{1.026}}
\multiput(158.097,36.498)(-.040810811,-.033441441){111}{\line(-1,0){.040810811}}
\multiput(153.567,31.99)(.0485,-.03341463){82}{\line(1,0){.0485}}
\put(158.048,9.269){\circle{1.026}}
\put(158.079,12.836){\line(0,-1){3.124}}
\put(153.318,13.368){\circle{1.026}}
\multiput(157.482,17.225)(-.042106796,-.033679612){103}{\line(-1,0){.042106796}}
\multiput(153.145,12.915)(.039669811,-.033707547){106}{\line(1,0){.039669811}}
\put(153.318,21.357){\circle{1.026}}
\multiput(157.614,25.319)(-.042097087,-.033679612){103}{\line(-1,0){.042097087}}
\multiput(153.147,21.009)(.038508621,-.033525862){116}{\line(1,0){.038508621}}
\put(55.434,17.362){\circle{1.026}}
\multiput(59.467,21.324)(-.040825243,-.033679612){103}{\line(-1,0){.040825243}}
\multiput(55.262,17.014)(.037212389,-.033486726){113}{\line(1,0){.037212389}}
\put(48.693,6.03){\makebox(0,0)[cc]{$B_{20}$}}
\put(59.927,6.03){\makebox(0,0)[cc]{$B_{21}$}}
\put(71.094,6.03){\makebox(0,0)[cc]{$B_{22}$}}
\put(82.263,6.03){\makebox(0,0)[cc]{$B_{23}$}}
\put(92.64,5.955){\makebox(0,0)[cc]{$B_{24}$}}
\put(104.908,6.023){\makebox(0,0)[cc]{$B_{25}$}}
\put(115.058,5.922){\makebox(0,0)[cc]{$B_{26}$}}
\put(127.42,5.971){\makebox(0,0)[cc]{$B_{27}$}}
\put(137.54,5.971){\makebox(0,0)[cc]{$B_{28}$}}
\put(147.42,5.971){\makebox(0,0)[cc]{$B_{29}$}}
\put(158.023,6.068){\makebox(0,0)[cc]{$B_{30}$}}
\put(133.047,25.124){\circle{1.026}}
\multiput(137.3,28.815)(-.04458,-.03373){100}{\line(-1,0){.04458}}
\multiput(133.083,24.671)(.041795918,-.03344898){98}{\line(1,0){.041795918}}
\put(133.168,15.099){\circle{1.026}}
\multiput(137.3,21.297)(-.033620155,-.046325581){129}{\line(0,-1){.046325581}}
\multiput(133.203,14.743)(.08361224,-.03344898){49}{\line(1,0){.08361224}}
\put(82.355,.25){\makebox(0,0)[cc]{Figure I}}
\put(82.417,9.408){\circle{1.026}}
\put(82.467,13.42){\circle{1.026}}
\put(82.417,17.432){\circle{1.026}}
\put(82.467,21.404){\circle{1.026}}
\put(82.448,20.999){\line(0,-1){3.124}}
\put(82.448,17.016){\line(0,-1){3.124}}
\put(82.448,12.975){\line(0,-1){3.124}}
\put(82.435,25.343){\circle{1.026}}
\put(82.417,24.938){\line(0,-1){3.124}}
\put(78.674,13.202){\circle{1.026}}
\multiput(78.537,12.852)(.034059406,-.033475248){101}{\line(1,0){.034059406}}
\put(82.515,29.084){\circle{1.026}}
\put(82.497,28.679){\line(0,-1){3.124}}
\put(82.484,33.023){\circle{1.026}}
\put(82.465,32.618){\line(0,-1){3.124}}
\put(78.033,28.537){\circle{1.026}}
\multiput(78.074,28.975)(.033478261,.035243478){115}{\line(0,1){.035243478}}
\multiput(78.074,28.164)(.04410588,-.03336471){85}{\line(1,0){.04410588}}
\multiput(81.933,17.364)(-.033510204,-.038612245){98}{\line(0,-1){.038612245}}
\put(78.033,17.396){\circle{1.026}}
\multiput(82.003,21.436)(-.036477064,-.033605505){109}{\line(-1,0){.036477064}}
\multiput(77.907,17.099)(.033694915,-.062898305){118}{\line(0,-1){.062898305}}
\put(127.577,9.302){\circle{1.026}}
\put(127.627,13.314){\circle{1.026}}
\put(127.577,17.326){\circle{1.026}}
\put(127.627,21.298){\circle{1.026}}
\put(127.608,20.893){\line(0,-1){3.124}}
\put(127.608,16.91){\line(0,-1){3.124}}
\put(127.608,12.869){\line(0,-1){3.124}}
\put(127.595,25.237){\circle{1.026}}
\put(127.577,24.832){\line(0,-1){3.124}}
\put(127.675,28.978){\circle{1.026}}
\put(127.657,28.573){\line(0,-1){3.124}}
\put(127.644,32.917){\circle{1.026}}
\put(127.625,32.512){\line(0,-1){3.124}}
\put(123.309,25.272){\circle{1.026}}
\put(123.193,28.431){\circle{1.026}}
\multiput(123.234,28.869)(.033478261,.035243478){115}{\line(0,1){.035243478}}
\multiput(123.184,28.025)(.0442381,-.03361905){84}{\line(1,0){.0442381}}
\put(123.1,17.329){\circle{1.026}}
\multiput(126.993,21.485)(-.03612037,-.033722222){108}{\line(-1,0){.03612037}}
\put(123.092,17.843){\line(0,1){0}}
\multiput(123.092,16.876)(.036801887,-.033650943){106}{\line(1,0){.036801887}}
\multiput(126.993,32.559)(-.033477477,-.060936937){111}{\line(0,-1){.060936937}}
\multiput(123.277,24.681)(.033477477,-.138612613){111}{\line(0,-1){.138612613}}
\put(105.069,9.408){\circle{1.026}}
\put(105.119,13.42){\circle{1.026}}
\put(105.069,17.432){\circle{1.026}}
\put(105.119,21.404){\circle{1.026}}
\put(105.1,20.999){\line(0,-1){3.124}}
\put(105.1,17.016){\line(0,-1){3.124}}
\put(105.1,12.975){\line(0,-1){3.124}}
\put(105.088,25.343){\circle{1.026}}
\put(105.069,24.938){\line(0,-1){3.124}}
\put(105.168,29.084){\circle{1.026}}
\put(105.149,28.679){\line(0,-1){3.124}}
\put(105.137,33.023){\circle{1.026}}
\put(105.118,32.618){\line(0,-1){3.124}}
\put(100.584,17.948){\line(0,1){0}}
\put(100.593,13.57){\circle{1.026}}
\put(100.5,24.682){\circle{1.026}}
\multiput(104.394,17.428)(-.037519231,-.033586538){104}{\line(-1,0){.037519231}}
\multiput(100.584,13.266)(.034539823,-.033548673){113}{\line(1,0){.034539823}}
\multiput(104.579,29.097)(-.0335,-.033508197){122}{\line(0,-1){.033508197}}
\multiput(100.399,24.34)(.04985366,-.03353659){82}{\line(1,0){.04985366}}
\put(100.593,20.482){\circle{1.026}}
\multiput(104.487,33.036)(-.033646552,-.103801724){116}{\line(0,-1){.103801724}}
\multiput(100.584,19.955)(.04759756,-.03353659){82}{\line(1,0){.04759756}}
\put(115.055,9.513){\circle{1.026}}
\put(115.105,13.525){\circle{1.026}}
\put(115.055,17.537){\circle{1.026}}
\put(115.105,21.509){\circle{1.026}}
\put(115.087,21.104){\line(0,-1){3.124}}
\put(115.087,17.121){\line(0,-1){3.124}}
\put(115.087,13.08){\line(0,-1){3.124}}
\put(115.074,25.448){\circle{1.026}}
\put(115.055,25.043){\line(0,-1){3.124}}
\put(115.154,29.189){\circle{1.026}}
\put(115.135,28.784){\line(0,-1){3.124}}
\put(115.123,33.128){\circle{1.026}}
\put(115.104,32.723){\line(0,-1){3.124}}
\put(110.57,18.053){\line(0,1){0}}
\put(110.579,13.675){\circle{1.026}}
\put(110.487,24.787){\circle{1.026}}
\multiput(114.379,17.533)(-.037509615,-.033596154){104}{\line(-1,0){.037509615}}
\multiput(110.57,13.37)(.034539823,-.033539823){113}{\line(1,0){.034539823}}
\multiput(114.565,29.202)(-.0335,-.033508197){122}{\line(0,-1){.033508197}}
\multiput(110.384,24.445)(.04986585,-.03353659){82}{\line(1,0){.04986585}}
\put(110.579,19.621){\circle{1.026}}
\multiput(115.029,32.695)(-.033711111,-.093592593){135}{\line(0,-1){.093592593}}
\multiput(110.57,19.242)(.033646552,-.083293103){116}{\line(0,-1){.083293103}}
\put(5.057,6.196){\makebox(0,0)[cc]{$B_{16}$}}
\put(14.812,6.121){\makebox(0,0)[cc]{$B_{17}$}}
\put(26.147,6.27){\makebox(0,0)[cc]{$B_{18}$}}
\put(37.667,6.121){\makebox(0,0)[cc]{$B_{19}$}}
\put(26.89,13.115){\circle{1.026}}
\put(26.94,17.127){\circle{1.026}}
\put(26.89,21.139){\circle{1.026}}
\put(26.94,25.111){\circle{1.026}}
\put(26.922,24.706){\line(0,-1){3.124}}
\put(26.922,20.723){\line(0,-1){3.124}}
\put(26.922,16.682){\line(0,-1){3.124}}
\put(26.909,29.05){\circle{1.026}}
\put(26.89,28.645){\line(0,-1){3.124}}
\put(23.148,25.085){\circle{1.026}}
\multiput(26.378,29.103)(-.033583333,-.037614583){96}{\line(0,-1){.037614583}}
\put(26.982,9.094){\circle{1.026}}
\put(26.869,12.723){\line(0,-1){3.258}}
\put(22.153,20.762){\circle{1.026}}
\multiput(26.448,24.728)(-.041038095,-.033666667){105}{\line(-1,0){.041038095}}
\put(22.153,17.049){\circle{1.026}}
\multiput(21.918,16.596)(.040810811,-.033441441){111}{\line(1,0){.040810811}}
\multiput(23.023,24.728)(.033535714,-.135741071){112}{\line(0,-1){.135741071}}
\multiput(22.028,20.397)(.04443299,-.033721649){97}{\line(1,0){.04443299}}
\multiput(26.558,24.816)(-.033732824,-.054648855){131}{\line(0,-1){.054648855}}
\put(38.307,13.046){\circle{1.026}}
\put(38.357,17.058){\circle{1.026}}
\put(38.307,21.07){\circle{1.026}}
\put(38.357,25.042){\circle{1.026}}
\put(38.338,24.637){\line(0,-1){3.124}}
\put(38.338,20.654){\line(0,-1){3.124}}
\put(38.338,16.613){\line(0,-1){3.124}}
\put(38.325,28.981){\circle{1.026}}
\put(38.307,28.576){\line(0,-1){3.124}}
\put(34.564,25.016){\circle{1.026}}
\multiput(37.794,29.034)(-.033583333,-.037614583){96}{\line(0,-1){.037614583}}
\put(38.398,9.025){\circle{1.026}}
\put(38.285,12.654){\line(0,-1){3.258}}
\put(33.569,20.693){\circle{1.026}}
\multiput(37.864,24.659)(-.041038095,-.033666667){105}{\line(-1,0){.041038095}}
\put(33.569,16.98){\circle{1.026}}
\multiput(33.334,16.527)(.040810811,-.033441441){111}{\line(1,0){.040810811}}
\multiput(34.439,24.659)(.033535714,-.135741071){112}{\line(0,-1){.135741071}}
\multiput(33.205,17.326)(.043981132,.033641509){106}{\line(1,0){.043981132}}
\multiput(33.51,20.406)(.033535211,-.049661972){142}{\line(0,-1){.049661972}}
\put(5.25,13.046){\circle{1.026}}
\put(5.3,17.058){\circle{1.026}}
\put(5.25,21.07){\circle{1.026}}
\put(5.3,25.042){\circle{1.026}}
\put(5.282,24.637){\line(0,-1){3.124}}
\put(5.282,20.654){\line(0,-1){3.124}}
\put(5.282,16.613){\line(0,-1){3.124}}
\put(5.269,28.981){\circle{1.026}}
\put(5.25,28.576){\line(0,-1){3.124}}
\put(1.508,25.016){\circle{1.026}}
\put(5.342,9.025){\circle{1.026}}
\put(5.229,12.654){\line(0,-1){3.258}}
\put(.513,16.98){\circle{1.026}}
\multiput(4.802,28.935)(-.03344,-.03485){100}{\line(0,-1){.03485}}
\multiput(4.7,21.316)(-.035818966,-.033534483){116}{\line(-1,0){.035818966}}
\multiput(.545,16.615)(.037672566,-.033707965){113}{\line(1,0){.037672566}}
\put(1.417,11.395){\circle{1.026}}
\multiput(1.357,11.023)(.05647541,-.03322951){61}{\line(1,0){.05647541}}
\multiput(4.802,24.802)(-.03344,-.1305){100}{\line(0,-1){.1305}}
\multiput(1.474,24.538)(.033510204,-.159816327){98}{\line(0,-1){.159816327}}
\put(14.793,13.115){\circle{1.026}}
\put(14.843,17.127){\circle{1.026}}
\put(14.793,21.139){\circle{1.026}}
\put(14.843,25.111){\circle{1.026}}
\put(14.824,24.706){\line(0,-1){3.124}}
\put(14.824,20.723){\line(0,-1){3.124}}
\put(14.824,16.682){\line(0,-1){3.124}}
\put(14.812,29.05){\circle{1.026}}
\put(14.793,28.645){\line(0,-1){3.124}}
\put(11.05,25.085){\circle{1.026}}
\multiput(14.28,29.103)(-.033572917,-.037614583){96}{\line(0,-1){.037614583}}
\put(14.884,9.094){\circle{1.026}}
\put(14.772,12.723){\line(0,-1){3.258}}
\put(10.055,20.762){\circle{1.026}}
\multiput(14.35,24.728)(-.041028571,-.033666667){105}{\line(-1,0){.041028571}}
\put(10.055,17.049){\circle{1.026}}
\multiput(9.82,16.596)(.040810811,-.033441441){111}{\line(1,0){.040810811}}
\multiput(10.925,24.728)(.033544643,-.135741071){112}{\line(0,-1){.135741071}}
\multiput(9.93,20.397)(.04443299,-.033721649){97}{\line(1,0){.04443299}}
\put(10.044,21.179){\line(0,1){0}}
\multiput(14.24,28.882)(-.03365625,-.088390625){128}{\line(0,-1){.088390625}}
\put(4.255,45.567){\circle{1.026}}
\put(4.305,49.579){\circle{1.026}}
\put(4.255,53.591){\circle{1.026}}
\put(4.305,57.563){\circle{1.026}}
\put(4.287,57.158){\line(0,-1){3.124}}
\put(4.287,53.175){\line(0,-1){3.124}}
\put(4.287,49.134){\line(0,-1){3.124}}
\put(4.274,61.502){\circle{1.026}}
\put(4.255,61.097){\line(0,-1){3.124}}
\put(.513,57.537){\circle{1.026}}
\multiput(3.743,61.555)(-.033583333,-.037614583){96}{\line(0,-1){.037614583}}
\put(.513,49.361){\circle{1.026}}
\multiput(3.743,53.379)(-.033583333,-.037614583){96}{\line(0,-1){.037614583}}
\multiput(.483,48.937)(.033608247,-.036340206){97}{\line(0,-1){.036340206}}
\put(.744,53.441){\circle{1.026}}
\multiput(.449,57.117)(.0335,-.0741){100}{\line(0,-1){.0741}}
\multiput(3.733,57.433)(-.033542553,-.038585106){94}{\line(0,-1){.038585106}}
\multiput(.533,53.058)(.033958333,-.033604167){96}{\line(1,0){.033958333}}
\put(14.847,45.407){\circle{1.026}}
\put(14.897,49.419){\circle{1.026}}
\put(14.847,53.431){\circle{1.026}}
\put(14.897,57.403){\circle{1.026}}
\put(14.878,56.998){\line(0,-1){3.124}}
\put(14.878,53.015){\line(0,-1){3.124}}
\put(14.878,48.974){\line(0,-1){3.124}}
\put(14.865,61.342){\circle{1.026}}
\put(14.847,60.937){\line(0,-1){3.124}}
\put(11.104,57.377){\circle{1.026}}
\multiput(14.334,61.395)(-.033583333,-.037614583){96}{\line(0,-1){.037614583}}
\put(11.104,49.201){\circle{1.026}}
\multiput(11.074,48.777)(.033608247,-.036340206){97}{\line(0,-1){.036340206}}
\put(11.335,53.281){\circle{1.026}}
\multiput(14.324,57.273)(-.033531915,-.038585106){94}{\line(0,-1){.038585106}}
\multiput(11.008,56.941)(.033715686,-.036519608){102}{\line(0,-1){.036519608}}
\multiput(11.152,52.872)(.033572917,-.035822917){96}{\line(0,-1){.035822917}}
\multiput(14.375,57.171)(-.03367,-.07566){100}{\line(0,-1){.07566}}
\put(26.274,45.503){\circle{1.026}}
\put(26.324,49.515){\circle{1.026}}
\put(26.274,53.527){\circle{1.026}}
\put(26.324,57.499){\circle{1.026}}
\put(26.305,57.094){\line(0,-1){3.124}}
\put(26.305,53.111){\line(0,-1){3.124}}
\put(26.305,49.07){\line(0,-1){3.124}}
\put(26.274,61.033){\line(0,-1){3.124}}
\put(22.532,57.473){\circle{1.026}}
\multiput(25.762,61.491)(-.033583333,-.037614583){96}{\line(0,-1){.037614583}}
\put(22.532,49.297){\circle{1.026}}
\multiput(25.762,53.315)(-.033583333,-.037614583){96}{\line(0,-1){.037614583}}
\put(22.763,53.377){\circle{1.026}}
\multiput(25.752,57.369)(-.033542553,-.038585106){94}{\line(0,-1){.038585106}}
\multiput(22.552,52.994)(.033958333,-.033604167){96}{\line(1,0){.033958333}}
\multiput(22.538,57.143)(.03360177,-.099911504){113}{\line(0,-1){.099911504}}
\multiput(22.467,48.948)(.042826087,-.033641304){92}{\line(1,0){.042826087}}
\put(37.378,45.503){\circle{1.026}}
\put(37.428,49.515){\circle{1.026}}
\put(37.378,53.527){\circle{1.026}}
\put(37.428,57.499){\circle{1.026}}
\put(37.409,57.094){\line(0,-1){3.124}}
\put(37.409,53.111){\line(0,-1){3.124}}
\put(37.409,49.07){\line(0,-1){3.124}}
\put(37.397,61.438){\circle{1.026}}
\put(37.378,61.033){\line(0,-1){3.124}}
\put(33.635,57.473){\circle{1.026}}
\put(33.635,49.297){\circle{1.026}}
\multiput(33.605,48.873)(.033608247,-.036340206){97}{\line(0,-1){.036340206}}
\put(33.867,53.377){\circle{1.026}}
\multiput(36.855,57.369)(-.033531915,-.038585106){94}{\line(0,-1){.038585106}}
\multiput(33.54,57.037)(.033715686,-.036519608){102}{\line(0,-1){.036519608}}
\multiput(33.683,52.968)(.033583333,-.035822917){96}{\line(0,-1){.035822917}}
\multiput(37.297,61.041)(-.033447619,-.107533333){105}{\line(0,-1){.107533333}}
\multiput(37.153,60.926)(-.039373626,-.033384615){91}{\line(-1,0){.039373626}}
\put(48.665,45.469){\circle{1.026}}
\put(48.715,49.481){\circle{1.026}}
\put(48.665,53.493){\circle{1.026}}
\put(48.715,57.465){\circle{1.026}}
\put(48.697,57.06){\line(0,-1){3.124}}
\put(48.697,53.077){\line(0,-1){3.124}}
\put(48.697,49.036){\line(0,-1){3.124}}
\put(48.684,61.404){\circle{1.026}}
\put(48.665,60.999){\line(0,-1){3.124}}
\put(44.923,57.439){\circle{1.026}}
\put(44.923,49.263){\circle{1.026}}
\put(45.154,53.343){\circle{1.026}}
\multiput(44.828,57.003)(.033715686,-.036519608){102}{\line(0,-1){.036519608}}
\multiput(48.44,60.892)(-.039362637,-.033384615){91}{\line(-1,0){.039362637}}
\multiput(48.082,53.556)(-.033642857,-.039183673){98}{\line(0,-1){.039183673}}
\multiput(44.785,48.913)(.034059406,-.033475248){101}{\line(1,0){.034059406}}
\multiput(48.153,57.453)(-.03338636,-.04168182){88}{\line(0,-1){.04168182}}
\multiput(45.073,52.982)(.033478261,-.037369565){92}{\line(0,-1){.037369565}}
\put(4.264,42.033){\makebox(0,0)[cc]{$B_{1}$}}
\put(14.709,42.115){\makebox(0,0)[cc]{$B_{2}$}}
\put(26.1,42.115){\makebox(0,0)[cc]{$B_{3}$}}
\put(37.419,42){\makebox(0,0)[cc]{$B_{4}$}}
\put(48.595,42.115){\makebox(0,0)[cc]{$B_{5}$}}
\put(59.843,42.057){\makebox(0,0)[cc]{$B_{6}$}}
\put(71.234,42.057){\makebox(0,0)[cc]{$B_{7}$}}
\put(82.625,42.057){\makebox(0,0)[cc]{$B_{8}$}}
\put(92.654,42.115){\makebox(0,0)[cc]{$B_{9}$}}
\put(104.369,42.235){\makebox(0,0)[cc]{$B_{10}$}}
\put(115.614,42.235){\makebox(0,0)[cc]{$B_{11}$}}
\put(126.962,42.154){\makebox(0,0)[cc]{$B_{12}$}}
\put(138.15,42.139){\makebox(0,0)[cc]{$B_{13}$}}
\put(6.76,45.082){\makebox(0,0)[cc]{$0$}}
\put(6.76,53.175){\makebox(0,0)[cc]{$b$}}
\put(6.892,49.076){\makebox(0,0)[cc]{$a$}}
\put(6.892,57.485){\makebox(0,0)[cc]{$c$}}
\put(6.892,61.374){\makebox(0,0)[cc]{$1$}}
\put(148.388,42.106){\makebox(0,0)[cc]{$B_{14}$}}
\put(157.717,42){\makebox(0,0)[cc]{$B_{15}$}}
\put(81.889,49.593){\circle{1.026}}
\put(81.939,53.605){\circle{1.026}}
\put(81.889,57.617){\circle{1.026}}
\put(81.939,61.589){\circle{1.026}}
\put(81.92,61.184){\line(0,-1){3.124}}
\put(81.92,57.201){\line(0,-1){3.124}}
\put(81.92,53.16){\line(0,-1){3.124}}
\put(81.908,65.528){\circle{1.026}}
\put(81.889,65.123){\line(0,-1){3.124}}
\put(78.147,61.563){\circle{1.026}}
\multiput(81.377,65.581)(-.033583333,-.037614583){96}{\line(0,-1){.037614583}}
\put(78.147,53.387){\circle{1.026}}
\multiput(81.377,57.405)(-.033583333,-.037614583){96}{\line(0,-1){.037614583}}
\multiput(78.117,52.963)(.033608247,-.036340206){97}{\line(0,-1){.036340206}}
\multiput(78.025,61.14)(.033565657,-.039242424){99}{\line(0,-1){.039242424}}
\put(81.98,45.572){\circle{1.026}}
\put(81.868,49.2){\line(0,-1){3.258}}
\put(78.164,48.55){\circle{1.026}}
\multiput(81.343,53.51)(-.033553191,-.047521277){94}{\line(0,-1){.047521277}}
\multiput(78.124,48.097)(.04437838,-.03337838){74}{\line(1,0){.04437838}}
\put(115.608,49.527){\circle{1.026}}
\put(115.658,53.539){\circle{1.026}}
\put(115.608,57.551){\circle{1.026}}
\put(115.658,61.523){\circle{1.026}}
\put(115.639,61.118){\line(0,-1){3.124}}
\put(115.639,57.135){\line(0,-1){3.124}}
\put(115.639,53.094){\line(0,-1){3.124}}
\put(115.627,65.462){\circle{1.026}}
\put(115.608,65.057){\line(0,-1){3.124}}
\put(111.865,61.497){\circle{1.026}}
\multiput(115.095,65.515)(-.033572917,-.037614583){96}{\line(0,-1){.037614583}}
\multiput(111.744,61.074)(.033565657,-.039242424){99}{\line(0,-1){.039242424}}
\put(115.699,45.506){\circle{1.026}}
\put(115.587,49.135){\line(0,-1){3.258}}
\put(111.883,48.484){\circle{1.026}}
\multiput(115.062,53.444)(-.033553191,-.047521277){94}{\line(0,-1){.047521277}}
\multiput(111.843,48.031)(.04437838,-.03337838){74}{\line(1,0){.04437838}}
\put(111.865,55.525){\circle{1.026}}
\put(93.125,49.527){\circle{1.026}}
\put(93.175,53.539){\circle{1.026}}
\put(93.125,57.551){\circle{1.026}}
\put(93.175,61.523){\circle{1.026}}
\put(93.157,61.118){\line(0,-1){3.124}}
\put(93.157,57.135){\line(0,-1){3.124}}
\put(93.157,53.094){\line(0,-1){3.124}}
\put(93.144,65.462){\circle{1.026}}
\put(93.125,65.057){\line(0,-1){3.124}}
\put(89.383,61.497){\circle{1.026}}
\multiput(92.613,65.515)(-.033583333,-.037614583){96}{\line(0,-1){.037614583}}
\multiput(89.262,61.074)(.033555556,-.039242424){99}{\line(0,-1){.039242424}}
\put(93.217,45.506){\circle{1.026}}
\put(93.104,49.135){\line(0,-1){3.258}}
\put(89.4,48.484){\circle{1.026}}
\multiput(92.579,53.444)(-.033553191,-.047521277){94}{\line(0,-1){.047521277}}
\multiput(89.36,48.031)(.04437838,-.03337838){74}{\line(1,0){.04437838}}
\put(89.383,55.525){\circle{1.026}}
\multiput(92.565,61.547)(-.033531915,-.058734043){94}{\line(0,-1){.058734043}}
\multiput(89.413,55.004)(.06259184,-.03338776){49}{\line(1,0){.06259184}}
\put(104.274,49.453){\circle{1.026}}
\put(104.324,53.465){\circle{1.026}}
\put(104.274,57.477){\circle{1.026}}
\put(104.324,61.449){\circle{1.026}}
\put(104.305,61.044){\line(0,-1){3.124}}
\put(104.305,57.061){\line(0,-1){3.124}}
\put(104.305,53.02){\line(0,-1){3.124}}
\put(104.293,65.388){\circle{1.026}}
\put(104.274,64.983){\line(0,-1){3.124}}
\put(100.532,61.423){\circle{1.026}}
\multiput(103.762,65.441)(-.033583333,-.037614583){96}{\line(0,-1){.037614583}}
\multiput(100.41,61)(.033565657,-.039242424){99}{\line(0,-1){.039242424}}
\put(104.365,45.432){\circle{1.026}}
\put(104.253,49.061){\line(0,-1){3.258}}
\put(100.549,48.41){\circle{1.026}}
\multiput(103.728,53.37)(-.033553191,-.047521277){94}{\line(0,-1){.047521277}}
\multiput(100.509,47.957)(.04437838,-.03337838){74}{\line(1,0){.04437838}}
\put(100.532,55.451){\circle{1.026}}
\multiput(103.714,65.29)(-.033531915,-.100787234){94}{\line(0,-1){.100787234}}
\multiput(100.477,55.066)(.03371875,-.05821875){96}{\line(0,-1){.05821875}}
\multiput(115.015,61.405)(-.033515464,-.055938144){97}{\line(0,-1){.055938144}}
\multiput(111.672,55.013)(.033696078,-.09327451){102}{\line(0,-1){.09327451}}
\put(126.739,49.63){\circle{1.026}}
\put(126.789,53.642){\circle{1.026}}
\put(126.739,57.654){\circle{1.026}}
\put(126.789,61.626){\circle{1.026}}
\put(126.77,61.221){\line(0,-1){3.124}}
\put(126.77,57.238){\line(0,-1){3.124}}
\put(126.77,53.197){\line(0,-1){3.124}}
\put(126.758,65.565){\circle{1.026}}
\put(126.739,65.16){\line(0,-1){3.124}}
\put(122.997,61.6){\circle{1.026}}
\multiput(126.227,65.618)(-.033583333,-.037614583){96}{\line(0,-1){.037614583}}
\put(126.83,45.609){\circle{1.026}}
\put(126.718,49.238){\line(0,-1){3.258}}
\put(123.014,48.587){\circle{1.026}}
\multiput(122.974,48.134)(.04437838,-.03337838){74}{\line(1,0){.04437838}}
\multiput(122.983,61.154)(.033726316,-.121884211){95}{\line(0,-1){.121884211}}
\put(121.892,53.564){\circle{1.026}}
\multiput(126.187,57.53)(-.041038095,-.033666667){105}{\line(-1,0){.041038095}}
\multiput(121.878,53.111)(.04199,-.03359){100}{\line(1,0){.04199}}
\multiput(126.187,61.685)(-.033578431,-.124784314){102}{\line(0,-1){.124784314}}
\put(149.442,49.523){\circle{1.026}}
\put(149.492,53.535){\circle{1.026}}
\put(149.442,57.547){\circle{1.026}}
\put(149.492,61.519){\circle{1.026}}
\put(149.473,61.114){\line(0,-1){3.124}}
\put(149.473,57.131){\line(0,-1){3.124}}
\put(149.473,53.09){\line(0,-1){3.124}}
\put(149.46,65.458){\circle{1.026}}
\put(149.442,65.053){\line(0,-1){3.124}}
\put(145.699,61.493){\circle{1.026}}
\multiput(148.929,65.511)(-.033583333,-.037614583){96}{\line(0,-1){.037614583}}
\multiput(145.578,61.07)(.033555556,-.039242424){99}{\line(0,-1){.039242424}}
\put(149.533,45.502){\circle{1.026}}
\put(149.42,49.131){\line(0,-1){3.258}}
\put(145.699,55.521){\circle{1.026}}
\multiput(148.849,61.401)(-.033515464,-.055938144){97}{\line(0,-1){.055938144}}
\multiput(145.504,55.009)(.033705882,-.09327451){102}{\line(0,-1){.09327451}}
\put(145.83,50.652){\circle{1.026}}
\multiput(148.79,57.337)(-.03346667,-.06961111){90}{\line(0,-1){.06961111}}
\multiput(145.658,50.301)(.1122069,-.0332414){29}{\line(1,0){.1122069}}
\put(139.003,49.5){\circle{1.026}}
\put(139.053,53.512){\circle{1.026}}
\put(139.003,57.524){\circle{1.026}}
\put(139.053,61.496){\circle{1.026}}
\put(139.034,61.091){\line(0,-1){3.124}}
\put(139.034,57.108){\line(0,-1){3.124}}
\put(139.034,53.067){\line(0,-1){3.124}}
\put(139.022,65.435){\circle{1.026}}
\put(139.003,65.03){\line(0,-1){3.124}}
\put(135.26,61.47){\circle{1.026}}
\multiput(138.49,65.488)(-.033572917,-.037614583){96}{\line(0,-1){.037614583}}
\put(139.094,45.479){\circle{1.026}}
\put(138.982,49.108){\line(0,-1){3.258}}
\put(135.278,48.457){\circle{1.026}}
\multiput(135.238,48.004)(.04437838,-.03337838){74}{\line(1,0){.04437838}}
\multiput(135.247,61.024)(.033715789,-.121884211){95}{\line(0,-1){.121884211}}
\multiput(138.45,61.555)(-.033578431,-.124784314){102}{\line(0,-1){.124784314}}
\put(60.662,45.453){\circle{1.026}}
\put(60.712,49.465){\circle{1.026}}
\put(60.662,53.477){\circle{1.026}}
\put(60.712,57.449){\circle{1.026}}
\put(60.693,57.044){\line(0,-1){3.124}}
\put(60.693,53.061){\line(0,-1){3.124}}
\put(60.693,49.02){\line(0,-1){3.124}}
\put(60.68,61.388){\circle{1.026}}
\put(60.662,60.983){\line(0,-1){3.124}}
\put(56.919,57.423){\circle{1.026}}
\multiput(60.149,61.441)(-.033583333,-.037614583){96}{\line(0,-1){.037614583}}
\put(56.919,49.247){\circle{1.026}}
\multiput(56.889,48.823)(.033608247,-.036340206){97}{\line(0,-1){.036340206}}
\put(57.15,53.327){\circle{1.026}}
\multiput(56.824,56.987)(.033715686,-.036519608){102}{\line(0,-1){.036519608}}
\multiput(56.967,52.918)(.033572917,-.035822917){96}{\line(0,-1){.035822917}}
\multiput(60.19,57.216)(-.03366,-.07565){100}{\line(0,-1){.07565}}
\multiput(60.069,61.278)(-.03339326,-.08268539){89}{\line(0,-1){.08268539}}
\put(70.662,45.203){\circle{1.026}}
\put(70.712,49.215){\circle{1.026}}
\put(70.662,53.227){\circle{1.026}}
\put(70.712,57.199){\circle{1.026}}
\put(70.693,56.794){\line(0,-1){3.124}}
\put(70.693,52.811){\line(0,-1){3.124}}
\put(70.693,48.77){\line(0,-1){3.124}}
\put(70.68,61.138){\circle{1.026}}
\put(70.662,60.733){\line(0,-1){3.124}}
\put(66.919,57.173){\circle{1.026}}
\multiput(70.149,61.191)(-.033583333,-.037614583){96}{\line(0,-1){.037614583}}
\put(66.919,48.997){\circle{1.026}}
\multiput(70.149,53.015)(-.033583333,-.037614583){96}{\line(0,-1){.037614583}}
\multiput(66.889,48.573)(.033608247,-.036340206){97}{\line(0,-1){.036340206}}
\put(67.15,53.077){\circle{1.026}}
\multiput(66.855,56.753)(.0335,-.0741){100}{\line(0,-1){.0741}}
\multiput(70.139,57.069)(-.033531915,-.038585106){94}{\line(0,-1){.038585106}}
\multiput(66.998,52.659)(.033478261,-.082228261){92}{\line(0,-1){.082228261}}
\put(135.044,55.828){\circle{1.026}}
\multiput(138.384,61.25)(-.033696078,-.049019608){102}{\line(0,-1){.049019608}}
\multiput(134.79,55.5)(.06201587,-.03373016){63}{\line(1,0){.06201587}}
\put(158.844,49.224){\circle{1.026}}
\put(158.894,53.236){\circle{1.026}}
\put(158.844,57.248){\circle{1.026}}
\put(158.894,61.22){\circle{1.026}}
\put(158.875,60.815){\line(0,-1){3.124}}
\put(158.875,56.832){\line(0,-1){3.124}}
\put(158.875,52.791){\line(0,-1){3.124}}
\put(158.863,65.159){\circle{1.026}}
\put(158.844,64.754){\line(0,-1){3.124}}
\put(155.102,61.194){\circle{1.026}}
\multiput(158.332,65.212)(-.033583333,-.037614583){96}{\line(0,-1){.037614583}}
\put(158.935,45.203){\circle{1.026}}
\put(158.823,48.832){\line(0,-1){3.258}}
\put(155.119,48.181){\circle{1.026}}
\multiput(155.079,47.728)(.04437838,-.03337838){74}{\line(1,0){.04437838}}
\put(154.107,56.871){\circle{1.026}}
\multiput(158.402,60.837)(-.041038095,-.033666667){105}{\line(-1,0){.041038095}}
\multiput(154.093,56.418)(.04199,-.03359){100}{\line(1,0){.04199}}
\multiput(154.755,48.639)(.033685714,.042933333){105}{\line(0,1){.042933333}}
\multiput(154.977,60.925)(.033578431,-.116117647){102}{\line(0,-1){.116117647}}
\put(26.367,61.474){\circle{1.026}}
\end{picture}

\end{center}

Let $B^{*}$ denote the dual of a lattice $B$. According to Bhavale \cite{bib1}, there are exactly thirty basic blocks (see Figure I) containing five comparable reducible elements, and having nullity three. Therefore we have the following result.
\begin{proposition}\label{bb4} If $B$ is the basic block associated to $\textbf{B}\in \mathscr{B}(j;5,3)$ where $j\geq 8$ then $B\in \{B_1,B_2,B_3,\ldots,B_{30}\}$\textnormal{(see Figure I)}.
\end{proposition}

Now out of thirty basic blocks there are seven of height four, namely $B_1$ to $B_{7}$, twelve of height five, namely $B_8$ to $B_{19}$, nine of height six, namely $B_{20}$ to $B_{28}$, and two of height seven, namely $B_{29}$, $B_{30}$. Also for $n\geq 8$, if $L \in \mathscr{L}(n;5,3,h)$ then $4 \leq h \leq 7$. For $1\leq i\leq 30$, let $\mathbb{B}_{i}(j;5,3)$ be the subclass of $\mathscr{B}(j;5,3)$ containing all maximal blocks of type $\textbf{B}\in\mathscr{B}(j;5,3)$ such that $B_i$(see Figure I) is the basic block associated to $\textbf{B}$.
\begin{remark}\label{h326}
\small{(1) For $j\geq 8$, $\mathscr{B}(j;5,3,4)=\displaystyle\dot\cup_{i=1}^{7}\mathbb{B}_i(j;5,3)$.
(2) For $j\geq 9$, $\mathscr{B}(j;5,3,5)=\displaystyle\dot\cup_{i=8}^{19}\mathbb{B}_i(j;5,3)$.
(3) For $j\geq 10$, $\mathscr{B}(j;5,3,6)=\displaystyle\dot\cup_{i=20}^{28}\mathbb{B}_i(j;5,3)$.
(4) For $j\geq 11$, $\mathscr{B}(j;5,3,7)=\displaystyle\dot\cup_{i=29}^{30}\mathbb{B}_i(j;5,3)$.
(5) For $j\geq8$, $\mathscr{B}(j;5,3)=\displaystyle\dot\cup_{h=4}^{7}\mathscr{B}(j;5,3,h)$.
(6) For $n\geq8$, $\mathscr{L}(n;5,3)=\displaystyle\dot\cup_{h=4}^{7}\mathscr{L}(n;5,3,h)$}.
\end{remark}

Now in the following, firstly we count the classes $\mathscr{B}(j;5,3,h)$ for $4 \leq h \leq 7$. Secondly, we count the classes $\mathscr{L}(n;5,3,h)$ for $4 \leq h \leq 7$, and thereby we count the class $\mathscr{L}(n;5,3)$.
\subsection{Counting of the class $\mathscr{B}(j;5,3,4)$}
 Here, we count the class $\mathscr{B}(j;5,3,4)$ by counting the classes $\mathbb{B}_i(j;5,3)$ for $i=1$ to $7$.
 \begin{proposition}\label{Bj5341}
For $j\geq 8$, $\displaystyle |\mathbb{B}_1(j;5,3)|=\displaystyle\sum_{s=1}^{j-7}\sum_{t=1}^{j-s-6}\binom{j-t-s-2}{4}$.
\end{proposition}
\begin{proof}
Let $\textbf{B}\in \mathbb{B}_1(j;5,3)$. Let $0<a<b<c<1$ be the reducible elements of $\textbf{B}$. As $B_1$(see Figure I) is the basic block associated to $\textbf{B}$, by Theorem \ref{redb}, $Red(B_1)=Red(\textbf{B})$ and $\eta(B_1)=\eta(\textbf{B})=3$. Observe that an adjunct representation of $B_1$ is given by $B_1=C]_{0}^{b}\{c_1\}]_{a}^{c}\{c_2\}]_{a}^{1}\{c_3\}$, where $C:0\prec a\prec b\prec c \prec 1$ is a $5$-chain. Also by Corollary \ref{maxchain}, $\textbf{B}$ has an adjunct representation $\textbf{B}=C_0]_{0}^{b}C_1]_{a}^{c}C_2]_{a}^{1}C_3$, where $C_0$ is a maximal chain containing all the reducible elements of $\textbf{B}$, and $C_1, C_2, C_3$ are chains. 

Observe that $\textbf{B}=(\textbf{B}'\oplus C')]_{a}^{1}C_3$, where $\textbf{B}'=\textbf{B}\cap[0,c]\in\mathscr{B}(i;4,2,3)$ for $i\geq 6$, $C'=\textbf{B}\cap(c,1]$. Let $|C_3|=s\geq 1$ and $|C'|=t\geq 1$ with $j=i+t+s\geq 8$. Suppose $\textbf{D}=(\textbf{D}'\oplus C'')]_{a}^{1}C_3'\in \mathbb{B}_1(j;5,3)$, where $\textbf{D}'=\textbf{D}\cap[0,c]\in\mathscr{B}(i;4,2,3)$ for $i\geq 6$, $C''=\textbf{D}\cap(c,1]$ with $|C''|=t\geq 1$, and $C_3'$ is a chain with $|C_3'|=s\geq 1$. Then $\textbf{B}\cong \textbf{D}$ if and only if $\textbf{B}'\cong \textbf{D}'$, $C'\cong C''$, and $C_3\cong C'_3$. To prove this, suppose $\textbf{B}\cong \textbf{D}$. As $|C_3|=|C'_3|=s$, $C_3\cong C'_3$, and hence $\textbf{B}\setminus C_3\cong \textbf{D}\setminus C'_3$, {\it{that is}}, $\textbf{B}'\oplus C'\cong \textbf{D}'\oplus C''$. Also $|C'|=|C''|=t$, $C'\cong C''$, and hence $(\textbf{B}'\oplus C')\setminus C'\cong (\textbf{D}'\oplus C'')\setminus C''$, {\it{that is}}, $\textbf{B}'\cong \textbf{D}'$. The converse is obvious.

 Now for fixed $t$ and $s$, there are $|\mathscr{B}(j-t-s;4,2,3)|$ maximal blocks in $\mathbb{B}_1(j;5,3)$ up to isomorphism. Therefore for fixed $s\geq 1$, $1\leq t=j-s-i\leq j-s-6$, since $i\geq 6$, and there are $\displaystyle\sum_{t=1}^{j-s-6}|\mathscr{B}(j-t-s;4,2,3)|$ maximal blocks in $\mathbb{B}_1(j;5,3)$ up to isomorphism. Further, $1\leq s=j-i-t\leq j-7$, since $i\geq 6$, $t\geq 1$, and there are $\displaystyle\sum_{s=1}^{j-7}\sum_{t=1}^{j-s-6}|\mathscr{B}(j-t-s;4,2,3)|$ maximal blocks in $\mathbb{B}_1(j;5,3)$ up to isomorphism. According to Aware and Bhavale \cite{bib17}(see Proposition 3.1), $|\mathscr{B}(i;4,2,3)|=\binom{i-2}{4}$. Thus there are $\displaystyle\sum_{s=1}^{j-7}\sum_{t=1}^{j-s-6}\binom{j-t-s-2}{4}$ maximal blocks in $\mathbb{B}_1(j;5,3)$ up to isomorphism.
\end{proof}
Note that $\textbf{B}\in\mathbb{B}_2(j;5,3)$ if and only if $\textbf{B*}\in\mathbb{B}_1(j;5,3)$. Therefore by Proposition \ref{Bj5341}, we have the following result.
\begin{corollary}\label{Bj5342}
For $j\geq 8$, $\displaystyle|\mathbb{B}_2(j;5,3)|= \displaystyle\sum_{s=1}^{j-7}\sum_{t=1}^{j-s-6}\binom{j-t-s-2}{4}$.
\end{corollary}

\begin{proposition}\label{Bj5343}For $j\geq 8$, $\displaystyle|\mathbb{B}_3(j;5,3)|=\displaystyle\sum_{s=1}^{j-7}\sum_{t=1}^{j-s-6}\binom{j-t-s-2}{4}$.
\end{proposition}
\begin{proof}
Let $\textbf{B}\in \mathbb{B}_3(j;5,3)$. Let $0<a<b<c<1$ be the reducible elements of $\textbf{B}$. As $B_3$(see Figure I) is the basic block associated to $\textbf{B}$, by Theorem \ref{redb}, $Red(B_3)=Red(\textbf{B})$ and $\eta(B_3)=\eta(\textbf{B})=3$. Observe that an adjunct representation of $B_3$ is given by $B_3=C]_{0}^{b}\{c_1\}]_{a}^{c}\{c_2\}]_{0}^{1}\{c_3\}$, where $C:0\prec a\prec b\prec c \prec 1$ is a $5$-chain. Also by Corollary \ref{maxchain}, $\textbf{B}$ has an adjunct representation $\textbf{B}=C_0]_{0}^{b}C_1]_{a}^{c}C_2]_{0}^{1}C_3$, where $C_0$ is a maximal chain containing all the reducible elements of $\textbf{B}$, and $C_1, C_2, C_3$ are chains.

Observe that $\textbf{B}=(\textbf{B}'\oplus C')]_{0}^{1}C_3$, where $\textbf{B}'=\textbf{B}\cap[0,c]\in\mathscr{B}(i;4,2,3)$ for $i\geq 6$, $C'=\textbf{B}\cap(c,1]$. Let $|C_3|=s\geq 1$ and $|C'|=t\geq 1$ with $j=i+t+s\geq 8$. Suppose $\textbf{D}=(\textbf{D}'\oplus C'')]_{0}^{1}C_3'\in \mathbb{B}_3(j;5,3)$, where $\textbf{D}'=\textbf{D}\cap[0,c]\in\mathscr{B}(i;4,2,3)$ for $i\geq 6$, $C''=\textbf{D}\cap(c,1]$ with $|C''|=t\geq 1$, and $|C_3'|=s\geq 1$. Then $\textbf{B}\cong \textbf{D}$ if and only if $\textbf{B}'\cong \textbf{D}'$, $C'\cong C''$, and $C_3\cong C'_3$.

Now observe that an adjunct representation of $\textbf{B}\in \mathbb{B}_3(j;5,3)$ is similar to an adjunct representation of $\textbf{B}\in \mathbb{B}_1(j;5,3)$ with an exception that the adjunct pair corresponds to $C_3$ is $(0,1)$ instead of $(a,1)$. That means the chain $C_3$ is joined to $0$ instead of $a$, and that doesn't hamper the counting in this case. Therefore the procedure of counting of the class $\mathbb{B}_3(j;5,3)$ is same as that of the class $\mathbb{B}_1(j;5,3)$. Moreover $|\mathbb{B}_3(j;5,3)|=|\mathbb{B}_1(j;5,3)|$.
\end{proof}
Note that $\textbf{B}\in\mathbb{B}_4(j;5,3)$ if and only if $\textbf{B*}\in\mathbb{B}_3(j;5,3)$. Therefore by Proposition \ref{Bj5343}, we have the following result.
\begin{corollary}\label{Bj5344}For $j\geq 8$, $\displaystyle|\mathbb{B}_4(j;5,3)|=\displaystyle\sum_{s=1}^{j-7}\sum_{t=1}^{j-s-6}\binom{j-t-s-2}{4}$.
\end{corollary}

\begin{proposition}\label{Bj5345}For $j\geq 8$, $\displaystyle|\mathbb{B}_5(j;5,3)|=\displaystyle\sum_{s=1}^{j-7}\sum_{t=1}^{j-s-6}\binom{j-t-s-2}{4}$.
\end{proposition}
\begin{proof}
Let $\textbf{B}\in \mathbb{B}_5(j;5,3)$. Let $0<a<b<c<1$ be the reducible elements of $\textbf{B}$. As $B_5$(see Figure I) is the basic block associated to $\textbf{B}$, by Theorem \ref{redb}, $Red(B_5)=Red(\textbf{B})$ and $\eta(B_5)=\eta(\textbf{B})=3$. Observe that an adjunct representation of $B_5$ is given by $B_5=C]_{0}^{b}\{c_1\}]_{a}^{c}\{c_2\}]_{b}^{1}\{c_3\}$, where $C:0\prec a\prec b\prec c \prec 1$ is a $5$-chain. Also by Corollary \ref{maxchain}, $\textbf{B}$ has an adjunct representation $\textbf{B}=C_0]_{0}^{b}C_1]_{a}^{c}C_2]_{b}^{1}C_3$, where $C_0$ is a maximal chain containing all the reducible elements of $\textbf{B}$, and $C_1, C_2, C_3$ are chains. 

Now observe that an adjunct representation of $\textbf{B}\in \mathbb{B}_5(j;5,3)$ is similar to an adjunct representation of $\textbf{B}\in \mathbb{B}_1(j;5,3)$ with an exception that the adjunct pair corresponds to $C_3$ is $(b,1)$ instead of $(a,1)$. Therefore the procedure of counting of the class $\mathbb{B}_5(j;5,3)$ is same as that of the class $\mathbb{B}_1(j;5,3)$. Moreover $|\mathbb{B}_5(j;5,3)|=|\mathbb{B}_1(j;5,3)|$.
\end{proof}

\begin{proposition}\label{Bj5346}For $j\geq 8$, $\displaystyle|\mathbb{B}_6(j;5,3)|=\displaystyle\binom{j-2}{6}$.
\end{proposition}
\begin{proof}
Let $\textbf{B}\in \mathbb{B}_6(j;5,3)$. Let $0<a<b<c<1$ be the reducible elements of $\textbf{B}$. As $B_6$(see Figure I) is the basic block associated to $\textbf{B}$, by Theorem \ref{redb}, $Red(B_6)=Red(\textbf{B})$ and $\eta(B_6)=\eta(\textbf{B})=3$. Observe that an adjunct representation of $B_6$ is given by $B_6=C]_{0}^{c}\{c_1\}]_{a}^{1}\{c_2\}]_{b}^{1}\{c_3\}$, where $C:0\prec a\prec b\prec c \prec 1$ is a $5$-chain. Also by Corollary \ref{maxchain}, $\textbf{B}$ has an adjunct representation $\textbf{B}=C_0]_{0}^{c}C_1]_{a}^{1}C_2]_{b}^{1}C_3$, where $C_0$ is a maximal chain containing all the reducible elements of $\textbf{B}$, and $C_1, C_2, C_3$ are chains. Let $s=|(0,a)\cap \textbf{B}|\geq 0$, $m=|(a,b)\cap \textbf{B}|\geq 0$, $t=|(b,c)\cap \textbf{B}|\geq 0$, $l=|(c,1)\cap \textbf{B}|\geq 0$, $p=|C_1|\geq 1$, $q=|C_2|\geq 1$, $r=|C_3|\geq 1$. Then $j=s+m+t+l+p+q+r+5\geq 8$.

Let $S=\{(s,m,t,l,p,q,r)~|~ s+m+t+l+p+q+r=j-5,~\textnormal{where}~ s,m,t,l\in\mathbb{N}\cup\{0\}~\textnormal{and}~p,q,r \in \mathbb{N}\}$. Now given any $\textbf{B}\in\mathbb{B}_6(j;5,3)$, there exists unique $(s,m,t,l,p,q,r)\in S$ and vice-versa. Also for $\textbf{B},\textbf{B}'\in\mathbb{B}_6(j;5,3)$, suppose there exist $(s,m,t,l,p,q,r),(s',m',t',l',p',q',r')\in S$ respectively. Then $\textbf{B}\cong\textbf{B}'$ if and only if $(s,m,t,l,p,q,r)=(s',m',t',l',p',q',r')$. Therefore the set $S$ is equivalent to the class $\mathbb{B}_6(j;5,3)$, and hence $|\mathbb{B}_6(j;5,3)|=|S|$. If $T=\{(x_1,x_2,x_3,x_4,x_5,x_6,x_7)~|~ x_1+x_2+x_3+x_4+x_5+x_6+x_7=j-8,~\textnormal{where}~ x_i \in \mathbb{N}\cup\{0\}, 1\leq i\leq 7\}$ then $|S|=|T|$. But $|T|$ is the number of all unordered non-negative integer solutions of the equation $x_1+x_2+x_3+x_4+x_5+x_6+x_7=j-8$. Therefore $|T|=\binom{(j-8)+7-1}{j-8}=\binom{j-2}{6}$. Thus $|\mathbb{B}_6(j;5,3)|=|S|=|T|=\binom{j-2}{6}$.
\end{proof}
Note that $\textbf{B}\in\mathbb{B}_7(j;5,3)$ if and only if $\textbf{B*}\in\mathbb{B}_6(j;5,3)$. Therefore by Proposition \ref{Bj5346}, we have the following result.
\begin{corollary}\label{Bj5347}For $j\geq 8$, $\displaystyle|\mathbb{B}_7(j;5,3)|=\displaystyle\binom{j-2}{6}$.
\end{corollary}

\begin{theorem}\label{Bj534}
For $j\geq 8$, $|\mathscr{B}(j;5,3,4)|=\displaystyle 2\binom{j-2}{6}+\displaystyle\sum_{s=1}^{j-7}\sum_{t=1}^{j-s-6}5\binom{j-t-s-2}{4}$.
\end{theorem}
\begin{proof} By (1) of Remark \ref{h326}, for $j\geq 8$, $\mathscr{B}(j;5,3,4)=\displaystyle\dot\cup_{i=1}^{7}\mathbb{B}_i(j;5,3)$. Therefore $\displaystyle|\mathscr{B}(j;5,3,4)|=\displaystyle\sum_{i=1}^{7}|\mathbb{B}_i(j;5,3)|$. Hence the proof follows from Proposition \ref{Bj5341}, Corollary \ref{Bj5342}, Proposition \ref{Bj5343}, Corollary \ref{Bj5344}, Proposition \ref{Bj5345}, Proposition \ref{Bj5346}, and Corollary \ref{Bj5347}.
\end{proof}
\subsection{Counting of the class $\mathscr{B}(j;5,3,5)$}
In this subsection, we count the classes $\mathbb{B}_i(j;5,3)$ for $i=8$ to $19$; Consequently, we count the class $\mathscr{B}(j;5,3,5)$.
\begin{proposition}\label{Bj5358}For $j\geq 9$, $|\mathbb{B}_8(j;5,3)|=\displaystyle\sum_{p=6}^{j-3}\binom{p-2}{4}P_{j-p-1}^{2}$.
\end{proposition}
\begin{proof}
Let $\textbf{B}\in \mathbb{B}_8(j;5,3)$. Let $0<a<b<c<1$ be the reducible elements of $\textbf{B}$. As $B_8$(see Figure I) is the basic block associated to $\textbf{B}$, by Theorem \ref{redb}, $Red(B_8)=Red(\textbf{B})$ and $\eta(B_8)=\eta(\textbf{B})=3$. Observe that an adjunct representation of $B_8$ is given by $B_8=C]_{0}^{b}\{c_1\}]_{a}^{c}\{c_2\}]_{c}^{1}\{c_3\}$, where $C:0\prec a\prec b\prec c\prec x \prec 1$ is a $6$-chain. Also by Corollary \ref{maxchain}, $\textbf{B}$ has an adjunct representation $\textbf{B}=C_0]_{0}^{b}C_1]_{a}^{c}C_2]_{c}^{1}C_3$, where $C_0$ is a maximal chain containing all the reducible elements of $\textbf{B}$, and $C_1, C_2, C_3$ are chains.

Observe that $\textbf{B}=\textbf{B}'\circ \textbf{B}''$, where $\textbf{B}'=\textbf{B}\cap[0,c]\in \mathscr{B}(p;4,2,3)$ for $p\geq 6$ and $\textbf{B}''=\textbf{B}\cap[c,1]\in \mathscr{B}(q;2,1,2)$ for $q\geq 4$ with $j=p+q-1\geq 9$. Note that $6\leq p=j-q+1\leq j-3$, since $q\geq 4$. By Corollary \ref{circ}, $|\mathbb{B}_8(j;5,3)|=\displaystyle\sum_{p=6}^{j-3}(|\mathscr{B}(p;4,2,3)|\times|\mathscr{B}(j-p+1;2,1,2)|)$. According to Aware and Bhavale \cite{bib17}(see Proposition 3.1), $|\mathscr{B}(p;4,2,3)|=\binom{p-2}{4}$. Also according to Thakare et al. \cite{bib14}(see Theorem 3.3) $|\mathscr{B}(q;2,1,2)|=P_{q-2}^{2}$. Thus $|\mathbb{B}_8(j;5,3)|=\displaystyle\sum_{p=6}^{j-3}\bigg{(}\binom{p-2}{4}\times P_{j-p-1}^{2}\bigg{)}$.
\end{proof}

Note that $\textbf{B}\in\mathbb{B}_9(j;5,3)$ if and only if $\textbf{B*}\in\mathbb{B}_8(j;5,3)$. Therefore by Proposition \ref{Bj5358}, we have the following result.
\begin{corollary}\label{Bj5359}
For $j\geq 9$, $|\mathbb{B}_9(j;5,3)|=\displaystyle\sum_{p=6}^{j-3}\binom{p-2}{4}P_{j-p-1}^{2}$.
\end{corollary}

\begin{proposition}\label{Bj53510}For $j\geq 9$,\\ $|\mathbb{B}_{10}(j;5,3)|=\displaystyle\sum_{p=0}^{j-9}\sum_{q=0}^{j-p-9}\sum_{r=0}^{j-p-q-9}\sum_{s=1}^{j-p-q-r-8}\sum_{t=1}^{j-p-q-r-s-7}P_{j-p-q-r-s-t-5}^{2}$.
\end{proposition}
\begin{proof}
Let $\textbf{B}\in \mathbb{B}_{10}(j;5,3)$. Let $0<a<b<c<1$ be the reducible elements of $\textbf{B}$. As $B_{10}$(see Figure I) is the basic block associated to $\textbf{B}$, by Theorem \ref{redb}, $Red(B_{10})=Red(\textbf{B})$ and $\eta(B_{10})=\eta(\textbf{B})=3$. Observe that an adjunct representation of $B_{10}$ is given by $B_{10}=C]_{0}^{b}\{c_1\}]_{a}^{1}\{c_2\}]_{c}^{1}\{c_3\}$, where $C:0\prec a\prec b\prec c\prec x \prec 1$ is a $6$-chain. Also by Corollary \ref{maxchain}, $\textbf{B}$ has an adjunct representation $\textbf{B}=C_0]_{0}^{b}C_1]_{a}^{1}C_2]_{c}^{1}C_3$, where $C_0$ is a maximal chain containing all the reducible elements of $\textbf{B}$, and $C_1, C_2, C_3$ are chains.

 Let $l=|(0,a)\cap \textbf{B}|\geq 0$, $m=|(a,b)\cap \textbf{B}|\geq 0$, $t=|(b,c)\cap \textbf{B}|\geq 0$, $p=|((c,1)\setminus C_3)\cap \textbf{B}|\geq 1$, $q=|C_1|\geq 1$, $r=|C_2|\geq 1$, $s=|C_3|\geq 1$. Then $j=l+m+t+p+q+r+s+5\geq 9$. Let $S=\{(l,m,t,p,q,r,s)~|~ l+m+t+p+q+r+s=j-5,p\geq s,~{\textnormal{where}}~ l,m,t\in\mathbb{N}\cup\{0\} ~{\textnormal{and}}~p,q,r,s \in \mathbb{N}\}$. Now given any $\textbf{B}\in\mathbb{B}_{10}(j;5,3)$, there exists unique $(l,m,t,p,q,r,s)\in S$ and vice-versa. Also for $\textbf{B},\textbf{B}'\in\mathbb{B}_{10}(j;5,3)$, suppose there exist $(l,m,t,p,q,r,s),(l',m',t',p',q',r',s')\in S$ respectively. Then $\textbf{B}\cong\textbf{B}'$ if and only if $(l,m,t,p,q,r,s)=(l',m',t',p',q',r',s')$. Therefore the set $S$ is equivalent to the class $\mathbb{B}_{10}(j;5,3)$, and hence $|\mathbb{B}_{10}(j;5,3)|=|S|$. Now let $T=\{(x_1,x_2,x_3,x_4,x_5,x_6)~|~ x_1+x_2+x_3+x_4+x_5+x_6=j-5, x_6=p+s ~{\textnormal{with}}~p\geq s\geq 1,~{\textnormal{where}}~ x_1,x_2,x_3 \in \mathbb{N}\cup\{0\}, ~{\textnormal{and}}~ x_4,x_5 \in \mathbb{N}\}$. Then $|S|=|T|$.

Now for fixed $x_1,x_2,x_3,x_4,x_5$, $2\leq x_6=j-x_1-x_2-x_3-x_4-x_5-5$, and there are $P_{x_6}^{2}$ maximal blocks in $\mathbb{B}_{10}(j;5,3)$ up to isomorphism. Also for fixed $x_1,x_2,x_3,x_4$, $1\leq x_5=j-x_1-x_2-x_3-x_4-x_6-5\leq j-x_1-x_2-x_3-x_4-7$, since $x_6\geq 2$, and there are $\displaystyle\sum_{x_5=1}^{j-x_1-x_2-x_3-x_4-7}P_{j-x_1-x_2-x_3-x_4-x_5-5}^{2}$ maximal blocks in $\mathbb{B}_{10}(j;5,3)$ up to isomorphism. Again, for fixed $x_1,x_2,x_3$, $1\leq x_4=j-x_1-x_2-x_3-x_5-x_6-5\leq j-x_1-x_2-x_3-8$, since $x_5\geq 1,x_6\geq 2$, and there are $\displaystyle\sum_{x_4=1}^{j-x_1-x_2-x_3-8}\sum_{x_5=1}^{j-x_1-x_2-x_3-x_4-7}P_{j-x_1-x_2-x_3-x_4-x_5-5}^{2}$ maximal blocks in $\mathbb{B}_{10}(j;5,3)$ up to isomorphism. Also, for fixed $x_1,x_2$, $0\leq x_3=j-x_1-x_2-x_4-x_5-x_6-5\leq j-x_1-x_2-9$, since $x_4\geq 1,x_5\geq 1,x_6\geq 2$, and there are $\displaystyle\sum_{x_3=0}^{j-x_1-x_2-9}\sum_{x_4=1}^{j-x_1-x_2-x_3-8}\sum_{x_5=1}^{j-x_1-x_2-x_3-x_4-7}P_{j-x_1-x_2-x_3-x_4-x_5-5}^{2}$ maximal blocks in $\mathbb{B}_{10}(j;5,3)$ up to isomorphism. Again, for fixed $x_1$, $0\leq x_2=j-x_1-x_3-x_4-x_5-x_6-5\leq j-x_1-9$, since $x_3\geq 0,x_4\geq 1,x_5\geq 1,x_6\geq 2$, and there are \\$\displaystyle\sum_{x_2=0}^{j-x_1-9}\sum_{x_3=0}^{j-x_1-x_2-9}\sum_{x_4=1}^{j-x_1-x_2-x_3-8}\sum_{x_5=1}^{j-x_1-x_2-x_3-x_4-7}P_{j-x_1-x_2-x_3-x_4-x_5-5}^{2}$ maximal blocks in $\mathbb{B}_{10}(j;5,3)$ up to isomorphism. Further, $0\leq x_1=j-x_2-x_3-x_4-x_5-x_6-5\leq j-9$, since $x_2\geq 0,x_3\geq 0,x_4\geq 1,x_5\geq 1,x_6\geq 2$, and there are \\$\displaystyle\sum_{x_1=0}^{j-9}\sum_{x_2=0}^{j-x_1-9}\sum_{x_3=0}^{j-x_1-x_2-9}\sum_{x_4=1}^{j-x_1-x_2-x_3-8}\sum_{x_5=1}^{j-x_1-x_2-x_3-x_4-7}P_{j-x_1-x_2-x_3-x_4-x_5-5}^{2}$ maximal blocks in \\$\mathbb{B}_{10}(j;5,3)$ up to isomorphism. Thus $|\mathbb{B}_6(j;5,3)|=|S|=|T|=$\\$\displaystyle\sum_{x_1=0}^{j-9}\sum_{x_2=0}^{j-x_1-9}\sum_{x_3=0}^{j-x_1-x_2-9}\sum_{x_4=1}^{j-x_1-x_2-x_3-8}\sum_{x_5=1}^{j-x_1-x_2-x_3-x_4-7}P_{j-x_1-x_2-x_3-x_4-x_5-5}^{2}$.
\end{proof}
Note that $\textbf{B}\in\mathbb{B}_{11}(j;5,3)$ if and only if $\textbf{B*}\in\mathbb{B}_{10}(j;5,3)$. Therefore by Proposition \ref{Bj53510}, we have the following result.
\begin{corollary}\label{Bj53511}For $j\geq 9$,\\ $|\mathbb{B}_{11}(j;5,3)|=\displaystyle\sum_{p=0}^{j-9}\sum_{q=0}^{j-p-9}\sum_{r=0}^{j-p-q-9}\sum_{s=1}^{j-p-q-r-8}\sum_{t=1}^{j-p-q-r-s-7}P_{j-p-q-r-s-t-5}^{2}$.
\end{corollary}

\begin{proposition}\label{Bj53512}For $j\geq 9$, $|\mathbb{B}_{12}(j;5,3)|=\displaystyle\sum_{t=1}^{j-8}\sum_{k=1}^{j-t-7}\displaystyle\sum_{p=1}^{j-t-k-6}\sum_{l=2}^{j-t-k-p-4}(l-1)P^{2}_{j-t-k-p-l-2}$.
\end{proposition}
\begin{proof}
Let $\textbf{B}\in \mathbb{B}_{12}(j;5,3)$. Let $0<a<b<c<1$ be the reducible elements of $\textbf{B}$. As $B_{12}$(see Figure I) is the basic block associated to $\textbf{B}$, by Theorem \ref{redb}, $Red(B_{12})=Red(\textbf{B})$ and $\eta(B_{12})=\eta(\textbf{B})=3$. Observe that an adjunct representation of $B_{12}$ is given by $B_{12}=C]_{0}^{c}\{c_1\}]_{a}^{b}\{c_2\}]_{a}^{1}\{c_3\}$, where $C:0\prec a\prec x\prec b\prec c \prec 1$ is a $6$-chain. Also by Corollary \ref{maxchain}, $\textbf{B}$ has an adjunct representation $\textbf{B}=C_0]_{0}^{c}C_1]_{a}^{b}C_2]_{a}^{1}C_3$, where $C_0$ is a maximal chain containing all the reducible element of $\textbf{B}$, and $C_1, C_2, C_3$ are chains.

Now observe that $\textbf{B}=(\textbf{B}'\oplus C')]_{a}^{1}C_3$ where $\textbf{B}'=\textbf{B}\cap[0,c]\in\mathscr{B}(i;4,2,4)$ for $i\geq 7$, $C'=\textbf{B}\cap(c,1]$. Let $|C_3|=k\geq 1$ and $|C'|=t\geq 1$ with $j=i+t+k\geq 9$. Suppose $\textbf{D}=(\textbf{D}'\oplus C'')]_{a}^{1}C_3\in \mathbb{B}_{12}(j;5,3)$, where $\textbf{D}'=\textbf{B}\cap[0,c]\in\mathscr{B}(i;4,2,4)$ for $i\geq 7$, $C''=\textbf{D}\cap(c,1]$ with $|C''|=t\geq 1$, and $|C_3'|=k\geq 1$. Then $\textbf{B}\cong \textbf{D}$ if and only if $\textbf{B}'\cong \textbf{D}'$, $C'\cong C''$, and $C_3\cong C_3'$. Note that the procedure of counting of the class $\mathbb{B}_{12}(j;5,3)$ is same as the counting of the class $\mathbb{B}_1(j;5,3)$ with an exception that the block $\textbf{B}'\in\mathscr{B}(i;4,2,4)$ is used instead of the block $\textbf{B}'\in\mathscr{B}(i;4,2,3)$.
 
 Now for fixed $t$ and $k$, there are $|\mathscr{B}(j-t-k;4,2,4)|$ maximal blocks in $\mathbb{B}_{12}(j;5,3)$ up to isomorphism. Therefore for fixed $t$, $1\leq k=j-i-t\leq j-t-7$, since $i\geq 7$, and there are $\displaystyle\sum_{k=1}^{j-t-7}|\mathscr{B}(j-t-k;4,2,4)|$ maximal blocks in $\mathbb{B}_{12}(j;5,3)$ up to isomorphism. Further, $1\leq t=j-i-k\leq j-8$, since $k\geq 1$, $i\geq 7$, and there are $\displaystyle\sum_{t=1}^{j-8}\sum_{k=1}^{j-t-7}|\mathscr{B}(j-t-k;4,2,4)|$ maximal blocks in $\mathbb{B}_{12}(j;5,3)$ up to isomorphism. According to Aware and Bhavale \cite{bib17}(see Proposition 3.2), $|\mathscr{B}(i;4,2,4)|=\displaystyle\sum_{p=1}^{i-6}\sum_{l=2}^{i-p-4}(l-1)P^{2}_{i-p-l-2}$. Thus there are $\displaystyle\sum_{t=1}^{j-8}\sum_{k=1}^{j-t-7}\displaystyle\sum_{p=1}^{j-t-k-6}\sum_{l=2}^{j-t-k-p-4}(l-1)P^{2}_{j-t-k-p-l-2}$ maximal blocks in $\mathbb{B}_{12}(j;5,3)$ up to isomorphism.
\end{proof}

Note that $\textbf{B}\in\mathbb{B}_{13}(j;5,3)$ if and only if $\textbf{B*}\in\mathbb{B}_{12}(j;5,3)$. Therefore by Proposition \ref{Bj53512}, we have the following result.
\begin{corollary}\label{Bj53513}
For $j\geq 9$, $|\mathbb{B}_{13}(j;5,3)|=\displaystyle\sum_{t=1}^{j-8}\sum_{k=1}^{j-t-7}\displaystyle\sum_{p=1}^{j-t-k-6}\sum_{l=2}^{j-t-k-p-4}(l-1)P^{2}_{j-t-k-p-l-2}$.
\end{corollary}

\begin{proposition}\label{Bj53514}
For $j\geq 9$, $|\mathbb{B}_{14}(j;5,3)|=\displaystyle\sum_{t=1}^{j-8}\sum_{k=1}^{j-t-7}\displaystyle\sum_{p=1}^{j-t-k-6}\sum_{l=2}^{j-t-k-p-4}(l-1)P^{2}_{j-t-k-p-l-2}$.
\end{proposition}
\begin{proof}
Let $\textbf{B}\in \mathbb{B}_{14}(j;5,3)$. Let $0<a<b<c<1$ be the reducible elements of $\textbf{B}$. As $B_{14}$(see Figure I) is the basic block associated to $\textbf{B}$, by Theorem \ref{redb}, $Red(B_{14})=Red(\textbf{B})$ and $\eta(B_{14})=\eta(\textbf{B})=3$. Observe that an adjunct representation of $B_{14}$ is given by $B_{14}=C]_{0}^{c}\{c_1\}]_{a}^{b}\{c_2\}]_{b}^{1}\{c_3\}$, where $C:0\prec a\prec x\prec b\prec c \prec 1$ is a $6$-chain. Also by Corollary \ref{maxchain}, $\textbf{B}$ has an adjunct representation $\textbf{B}=C_0]_{0}^{c}C_1]_{a}^{b}C_2]_{b}^{1}C_3$, where $C_0$ is a maximal chain containing all the reducible element of $\textbf{B}$, and $C_1, C_2, C_3$ are chains.

Now observe that an adjunct representation of $\textbf{B}\in \mathbb{B}_{14}(j;5,3)$ is similar to an adjunct representation of $\textbf{B}\in \mathbb{B}_{12}(j;5,3)$(see proof of the Proposition \ref{Bj53512}) with an exception that the adjunct pair corresponds to $C_3$ is $(b,1)$ instead of $(a,1)$. Therefore the procedure of counting of the class $\mathbb{B}_{14}(j;5,3)$ is same as that of the class $\mathbb{B}_{12}(j;5,3)$. Moreover $|\mathbb{B}_{14}(j;5,3)|=|\mathbb{B}_{12}(j;5,3)|$.
\end{proof}

Note that $\textbf{B}\in\mathbb{B}_{15}(j;5,3)$ if and only if $\textbf{B*}\in\mathbb{B}_{14}(j;5,3)$. Therefore by Proposition \ref{Bj53514}, we have the following result.
\begin{corollary}\label{Bj53515}
For $j\geq 9$, $|\mathbb{B}_{15}(j;5,3)|=\displaystyle\sum_{t=1}^{j-8}\sum_{k=1}^{j-t-7}\displaystyle\sum_{p=1}^{j-t-k-6}\sum_{l=2}^{j-t-k-p-4}(l-1)P^{2}_{j-t-k-p-l-2}$.
\end{corollary}

\begin{proposition}\label{Bj53516}
For $j\geq 9$, $|\mathbb{B}_{16}(j;5,3)|=\displaystyle\sum_{t=1}^{j-8}\sum_{k=1}^{j-t-7}\displaystyle\sum_{p=1}^{j-t-k-6}\sum_{l=2}^{j-t-k-p-4}(l-1)P^{2}_{j-t-k-p-l-2}$.
\end{proposition}
\begin{proof}
Let $\textbf{B}\in \mathbb{B}_{16}(j;5,3)$. Let $0<a<b<c<1$ be the reducible elements of $\textbf{B}$. As $B_{16}$(see Figure I) is the basic block associated to $\textbf{B}$, by Theorem \ref{redb}, $Red(B_{16})=Red(\textbf{B})$ and $\eta(B_{16})=\eta(\textbf{B})=3$. Observe that an adjunct representation of $B_{16}$ is given by $B_{16}=C]_{0}^{c}\{c_1\}]_{a}^{b}\{c_2\}]_{0}^{1}\{c_3\}$, where $C:0\prec a\prec x\prec b\prec c \prec 1$ is a $6$-chain. Also by Corollary \ref{maxchain}, $\textbf{B}$ has an adjunct representation $\textbf{B}=C_0]_{0}^{c}C_1]_{a}^{b}C_2]_{0}^{1}C_3$, where $C_0$ is a maximal chain containing all the reducible element of $\textbf{B}$, and $C_1, C_2, C_3$ are chains.

Now observe that an adjunct representation of $\textbf{B}\in \mathbb{B}_{16}(j;5,3)$ is similar to an adjunct representation of $\textbf{B}\in \mathbb{B}_{14}(j;5,3)$(see proof of the Proposition \ref{Bj53514}) with an exception that the adjunct pair corresponds to $C_3$ is $(0,1)$ instead of $(b,1)$. Therefore the procedure of counting of the class $\mathbb{B}_{18}(j;5,3)$ is same as that of the class $\mathbb{B}_{14}(j;5,3)$. Moreover $|\mathbb{B}_{16}(j;5,3)|=|\mathbb{B}_{14}(j;5,3)|$.
\end{proof}
Note that $\textbf{B}\in\mathbb{B}_{17}(j;5,3)$ if and only if $\textbf{B*}\in\mathbb{B}_{16}(j;5,3)$. Therefore by Proposition \ref{Bj53516}, we have the following result.
\begin{corollary}\label{Bj53517}
For $j\geq 9$, $|\mathbb{B}_{17}(j;5,3)|=\displaystyle\sum_{t=1}^{j-8}\sum_{k=1}^{j-t-7}\displaystyle\sum_{p=1}^{j-t-k-6}\sum_{l=2}^{j-t-k-p-4}(l-1)P^{2}_{j-t-k-p-l-2}$.
\end{corollary}

\begin{proposition}\label{Bj53518}
For $j\geq 9$, $|\mathbb{B}_{18}(j;5,3)|=\displaystyle\sum_{t=1}^{j-8}\sum_{h=2}^{j-t-6}\displaystyle\sum_{l=1}^{j-h-t-5}\sum_{p=1}^{j-h-t-l-4}(h-1)P^{2}_{j-h-t-l-p-2}$.
\end{proposition}
\begin{proof}
Let $\textbf{B}\in \mathbb{B}_{18}(j;5,3)$. Let $0<a<b<c<1$ be the reducible elements of $\textbf{B}$. As $B_{18}$(see Figure I) is the basic block associated to $\textbf{B}$, by Theorem \ref{redb}, $Red(B_{18})=Red(\textbf{B})$ and $\eta(B_{18})=\eta(\textbf{B})=3$. Observe that an adjunct representation of $B_{18}$ is given by $B_{18}=C]_{a}^{c}\{c_1\}]_{b}^{c}\{c_2\}]_{0}^{1}\{c_3\}$, where $C:0\prec a\prec b\prec x\prec c \prec 1$ is a $6$-chain. Also by Corollary \ref{maxchain}, $\textbf{B}$ has an adjunct representation $\textbf{B}=C_0]_{a}^{c}C_1]_{b}^{c}C_2]_{0}^{1}C_3$, where $C_0$ is a maximal chain containing all the reducible element of $\textbf{B}$, and $C_1, C_2, C_3$ are chains.

Now observe that $\textbf{B}=(C'\oplus\textbf{B}'\oplus C'')]_{0}^{1}C_3$, where $\textbf{B}'=\textbf{B}\cap[a,c]\in\mathscr{B}(i;3,2,3)$ for $i\geq 6$, $C'=\textbf{B}\cap[0,a)$, $C''=\textbf{B}\cap(c,1]$. Let $|C'|=h_1\geq 1$, $|C''|=h_2\geq 1$, $|C_3|=t\geq 1$ and $h_1+h_2=h\geq 2$ with $j=i+h+t\geq 9$. Suppose $\textbf{D}=(E\oplus\textbf{D}'\oplus E')]_{0}^{1}C_3'\in \mathbb{B}_{18}(j;5,3)$, where $\textbf{D}'=\textbf{D}\cap[a,c]\in\mathscr{B}(i;3,2,3)$ for $i\geq 6$, $E=\textbf{D}\cap[0,a)$, $E'=\textbf{D}\cap(c,1]$ with $|E|=h_1\geq 1$,$|E'|=h_2\geq 1$, and $C_3'$ is a chain with $|C_3'|=t\geq 1$. Then $\textbf{B}\cong \textbf{D}$ if and only if $\textbf{B}'\cong \textbf{D}'$, $C'\cong E$, $C''\cong E'$, and $C_3\cong C'_3$.

 Note that $h-2$(excluding $0$ and $1$) elements can be distributed over the chains $C'$ and $C''$ in $h-2+1=h-1$ ways. Now for fixed $h_1,h_2$ and $t$, there are $(h-1)|\mathscr{B}(j-h-t;3,2,3)|$ maximal blocks in $\mathbb{B}_{18}(j;5,3)$ up to isomorphism. Therefore for fixed $t\geq 1$, $2\leq h=j-i-t\leq j-t-6$, since $i\geq 6$, and there are $\displaystyle\sum_{h=2}^{j-t-6}(h-1)|\mathscr{B}(j-h-t;3,2,3)|$ maximal blocks in $\mathbb{B}_{18}(j;5,3)$ up to isomorphism. Further, $1\leq t=j-h-i\leq j-8$, since $i\geq 6$, $h\geq 2$, and there are $\displaystyle\sum_{t=1}^{j-8}\sum_{h=2}^{j-t-6}(h-1)|\mathscr{B}(j-h-t;3,2,3)|$ maximal blocks in $\mathbb{B}_{18}(j;5,3)$ up to isomorphism. According to Bhavale and Aware \cite{bib2}(see Corollary 3.6 and put k=1), $|\mathscr{B}(i;3,2,3)|=\displaystyle\sum_{l=1}^{i-5}\sum_{p=1}^{i-l-4}P^{2}_{i-l-p-2}$. Therefore there are $\displaystyle\sum_{t=1}^{j-8}\sum_{h=2}^{j-t-6}(h-1)\bigg{(}\displaystyle\sum_{l=1}^{j-h-t-5}\sum_{p=1}^{j-h-t-l-4}P^{2}_{j-h-t-l-p-2}\bigg{)}$ maximal blocks in $\mathbb{B}_{18}(j;5,3)$ up to isomorphism.  Thus there are $\displaystyle\sum_{t=1}^{j-8}\sum_{h=2}^{j-t-6}\displaystyle\sum_{l=1}^{j-h-t-5}\sum_{p=1}^{j-h-t-l-4}(h-1)P^{2}_{j-h-t-l-p-2}$ maximal blocks in $\mathbb{B}_{18}(j;5,3)$ up to isomorphism.
\end{proof}
Note that $\textbf{B}\in\mathbb{B}_{19}(j;5,3)$ if and only if $\textbf{B*}\in\mathbb{B}_{18}(j;5,3)$. Therefore by Proposition \ref{Bj53518}, we have the following result.
\begin{corollary}\label{Bj53519}
For $j\geq 9$, $|\mathbb{B}_{19}(j;5,3)|=\displaystyle\sum_{t=1}^{j-8}\sum_{h=2}^{j-t-6}\displaystyle\sum_{l=1}^{j-h-t-5}\sum_{p=1}^{j-h-t-l-4}(h-1)P^{2}_{j-h-t-l-p-2}$.
\end{corollary}

\begin{theorem}\label{Bj535}
For $j\geq 9$,\\ \scriptsize{$|\mathscr{B}(j;5,3,5)|=\displaystyle\sum_{p=6}^{j-3}2\binom{p-2}{4}P_{j-p-1}^{2}+\displaystyle\sum_{p=0}^{j-9}\sum_{q=0}^{j-p-9}\sum_{r=0}^{j-p-q-9}\sum_{s=1}^{j-p-q-r-8}\sum_{t=1}^{j-p-q-r-s-7}2P_{j-p-q-r-s-t-5}^{2}+$\\$\displaystyle\sum_{t=1}^{j-8}\sum_{k=1}^{j-t-7}\sum_{p=1}^{j-t-k-6}\sum_{l=2}^{j-t-k-p-4}6(l-1)P^{2}_{j-t-k-p-l-2}+\displaystyle\sum_{t=1}^{j-8}\sum_{h=2}^{j-t-6}\sum_{l=1}^{j-h-t-5}\sum_{p=1}^{j-h-t-l-4}2(h-1)P^{2}_{j-h-t-l-p-2}$}.
\end{theorem}
\begin{proof}By (2) of Remark \ref{h326}, For $j\geq 9$, $\mathscr{B}(j;5,3,5)=\displaystyle\dot\cup_{i=8}^{19}\mathbb{B}_i(j;5,5)$. Therefore $\displaystyle|\mathscr{B}(j;5,3,5)|=\displaystyle\sum_{i=8}^{19}|\mathbb{B}_i(j;5,3)|$. Therefore the proof follows from Proposition \ref{Bj5358}, Corollary \ref{Bj5359}, Proposition \ref{Bj53510}, Corollary \ref{Bj53511}, Proposition \ref{Bj53512}, Corollary \ref{Bj53513}, Proposition \ref{Bj53514}, Proposition \ref{Bj53515}, Corollary \ref{Bj53516}, Proposition \ref{Bj53517}, Corollary \ref{Bj53518}, and Proposition \ref{Bj53519}.
\end{proof}
\subsection{Counting of the class $\mathbb{B}(j;5,3,6)$}
In this subsection, we count the classes $\mathbb{B}_i(j;5,3)$ for $i=20$ to $28$; Consequently, we count the class $\mathscr{B}(j;5,3,6)$. Let $\mathscr{L}'(n;r,k,h)$ be the subclass of $\mathscr{L}(n;r,k,h)$ such that every $L\in \mathscr{L}'(n;r,k,h)$ is of the type $L=\textbf{B}\oplus C$, where $\textbf{B}\in\mathscr{B}(m;r,k,h)$ and $C$ is a chain containing $n-m$ elements.
\begin{proposition}\label{Bj53520}
For $j\geq 10$, $\displaystyle|\mathbb{B}_{20}(j;5,3)|=\displaystyle\sum_{u=4}^{j-6}\sum_{r=0}^{u-4}\sum_{l=1}^{j-u-5}\sum_{s=1}^{j-u-l-4}P_{u-r-2}^{2}P^{2}_{j-u-l-s-2}$.
\end{proposition}
\begin{proof}
Let $\textbf{B}\in \mathbb{B}_{20}(j;5,3)$. Let $0<a<b<c<1$ be the reducible elements of $\textbf{B}$. As $B_{20}$(see Figure I) is the basic block associated to $\textbf{B}$, by Theorem \ref{redb}, $Red(B_{20})=Red(\textbf{B})$ and $\eta(B_{20})=\eta(\textbf{B})=3$. Observe that an adjunct representation of $B_{20}$ is given by $B_{20}=C]_{0}^{a}\{c_1\}]_{b}^{c}\{c_2\}]_{b}^{1}\{c_3\}$, where $C:0\prec x\prec a\prec b\prec y\prec c\prec 1$ is a $7$-chain. Also by Corollary \ref{maxchain}, $\textbf{B}$ has an adjunct representation $\textbf{B}=C_0]_{0}^{a}C_1]_{b}^{c}C_2]_{b}^{1}C_3$, where $C_0$ is a maximal chain containing all the reducible elements of $\textbf{B}$, and $C_1, C_2, C_3$ are chains. 

Now observe that $\textbf{B}=\textbf{B}'\oplus C'\oplus\textbf{B}''$, where $\textbf{B}'=[0,a]\cap \textbf{B}\in \mathscr{B}(p;2,1,2)$ with $p\geq 4$, $ \textbf{B}''=[b,1]\cap \textbf{B} \in \mathscr{B}(q;3,2,3)$ with $q\geq 6$, $C'$ is a chain with $|C'|=r\geq 0$, and $j=p+q+r\geq 10$. Suppose $\textbf{D}=\textbf{D}'\oplus C''\oplus\textbf{D}''\in\mathbb{B}_{20}(j;5,3)$, where $\textbf{D}'=[0,a]\cap \textbf{B}\in \mathscr{B}(p;2,1,2)$ with $p\geq 4$, $ \textbf{D}''=[b,1]\cap \textbf{B}\in \mathscr{B}(q;3,2,3)$ with $q\geq 6$, $C''$ is a chain with $|C''|=r\geq 0$. Then $\textbf{B}\cong \textbf{D}$ if and only if $\textbf{B}'\cong \textbf{D}'$, $\textbf{B}''\cong \textbf{D}''$, and $C'\cong C''$. 

Let $L=\textbf{B}'\oplus C'\in \mathscr{L}'(u;2,1,2)$, where $u=p+r\geq 4$. Then $\textbf{B}=L\oplus\textbf{B}''$. Now for fixed $r$, there are $|\mathscr{B}(u-r;2,1,2)|$ lattices of type $L$ in $\mathscr{L}'(u;2,1,2)$ up to isomorphism. Further, $0\leq r=u-p\leq u-4$, since $p\geq 4$, and there are $\displaystyle\sum_{r=0}^{u-4}|\mathscr{B}(u-r;2,1)|$ lattices of type $L$ in $\mathscr{L}'(u;2,1,2)$ up to isomorphism. According to Thakare et al. \cite{bib14}(see Theorem 3.3) $|\mathscr{B}(p;2,1,2)|=P_{p-2}^{2}$. Therefore there are $\displaystyle\sum_{r=0}^{u-4}P_{u-r-2}^{2}$ lattices of type $L$ in $\mathscr{L}'(u;2,1,2)$ up to isomorphism. 

Now $\textbf{B}=L\oplus\textbf{B}''$ and $4\leq u=j-q\leq j-6$, since $q\geq 6$. Therefore by Lemma \ref{oplus}, there are $\displaystyle\sum_{u=4}^{j-6}(|\mathscr{L}'(u;2,1,2)|\times|\mathscr{B}(j-u;3,2,3)|)$ maximal blocks in $\mathbb{B}_{20}(j;5,3)$ up to isomorphism.  Also according to Bhavale and Aware \cite{bib2}(see Corollary 3.6 and put k=1), $|\mathscr{B}(q;3,2,3)|=\displaystyle\sum_{l=1}^{q-5}\sum_{i=1}^{q-l-4}P^{2}_{q-l-i-2}$. Therefore there are $\displaystyle\sum_{u=4}^{j-6}\bigg{(}\displaystyle\sum_{r=0}^{u-4}P_{u-r-2}^{2}\times\displaystyle\sum_{l=1}^{j-u-5}\sum_{i=1}^{j-u-l-4}P^{2}_{j-u-l-i-2}\bigg{)}$ maximal blocks in $\mathbb{B}_{20}(j;5,3)$ up to isomorphism. Thus there are $\displaystyle\sum_{u=4}^{j-6}\displaystyle\sum_{r=0}^{u-4}\sum_{l=1}^{j-u-5}\sum_{i=1}^{j-u-l-4}P_{u-r-2}^{2}P^{2}_{j-u-l-i-2}$ maximal blocks in $\mathbb{B}_{20}(j;5,3)$ up to isomorphism.
\end{proof}

Note that $\textbf{B}\in\mathbb{B}_{21}(j;5,3)$ if and only if $\textbf{B*}\in\mathbb{B}_{20}(j;5,3)$. Therefore using Proposition \ref{Bj53520}, we have the following result.
\begin{corollary}\label{Bj53521}
For $j\geq 10$, $\displaystyle|\mathbb{B}_{21}(j;5,3)|=\displaystyle\sum_{u=4}^{j-6}\displaystyle\sum_{r=0}^{u-4}\sum_{l=1}^{j-u-5}\sum_{i=1}^{j-u-l-4}P_{u-r-2}^{2}P^{2}_{j-u-l-i-2}$.
\end{corollary}

\begin{proposition}\label{Bj53522}
For $j\geq 10$, $\displaystyle|\mathbb{B}_{22}(j;5,3)|=\displaystyle\sum_{u=4}^{j-6}\displaystyle\sum_{r=0}^{u-4}\sum_{l=1}^{j-u-5}\sum_{i=1}^{j-u-l-4}P_{u-r-2}^{2}P^{2}_{j-u-l-i-2}$.
\end{proposition}
\begin{proof}
Let $\textbf{B}\in \mathbb{B}_{22}(j;5,3)$. Let $0<a<b<c<1$ be the reducible elements of $\textbf{B}$. As $B_{22}$(see Figure I) is the basic block associated to $\textbf{B}$, by Theorem \ref{redb}, $Red(B_{22})=Red(\textbf{B})$ and $\eta(B_{22})=\eta(\textbf{B})=3$. Observe that an adjunct representation of $B_{22}$ is given by $B_{22}=C]_{0}^{a}\{c_1\}]_{b}^{1}\{c_2\}]_{c}^{1}\{c_3\}$, where $C:0\prec x\prec a\prec b\prec c\prec y\prec 1$ is a $7$-chain. Also by Corollary \ref{maxchain}, $\textbf{B}$ has an adjunct representation $\textbf{B}=C_0]_{0}^{a}C_1]_{b}^{1}C_2]_{c}^{1}C_3$, where $C_0$ is a maximal chain containing all the reducible elements of $\textbf{B}$, and $C_1, C_2, C_3$ are chains. 

$\textbf{B}=$ where $\textbf{B}'\in $ with $p\geq 4$, $ \textbf{B}''\in $ with $q\geq 6$, $C'$ is a chain with $|C'|=r\geq 0$, and $j=p+q+r\geq 10$.

Now observe that for $\textbf{B}\in \mathbb{B}_{22}(j;5,3)$, $\textbf{B}=\textbf{B}_1\oplus C'\oplus\textbf{B}_2$, where $\textbf{B}_1=\textbf{B}\cap[0,a]\in \mathscr{B}(p;2,1,2)$ with $p\geq 4$, $ \textbf{B}_2=\textbf{B}\cap[b,1]\in \mathscr{B}(q;3,2,3)$ with $q\geq 6$, $C'$ is a chain with $|C'|=r\geq 0$, and $j=p+q+r\geq 10$. Also for $\textbf{B}'\in \mathbb{B}_{20}(j;5,3)$, $\textbf{B}'=\textbf{B}_1'\oplus C''\oplus\textbf{B}_2'$, where $\textbf{B}_1'=\textbf{B}'\cap[0,a]\in \mathscr{B}(p;2,1,2)$ with $p\geq 4$, $ \textbf{B}_2'=\textbf{B}'\cap[b,1]\in \mathscr{B}(q;3,2,3)$ with $q\geq 6$, $C''$ is a chain with $|C''|=r\geq 0$, and $j=p+q+r\geq 10$. Note that $\textbf{B}_1\cong \textbf{B}_1'$ and $\textbf{B}_2\cong (\textbf{B}_2')^*$. Therefore counting of the class $\mathbb{B}_{22}(j;5,3)$ is same as that of the class $\mathbb{B}_{20}(j;5,3)$. Thus $|\mathbb{B}_{22}(j;5,3)|=|\mathbb{B}_{20}(j;5,3)|$. 
\end{proof}
Note that $\textbf{B}\in\mathbb{B}_{23}(j;5,3)$ if and only if $\textbf{B*}\in\mathbb{B}_{22}(j;5,3)$. Therefore using Proposition \ref{Bj53522}, we have the following result.
\begin{corollary}\label{Bj53523}
For $j\geq 10$, $\displaystyle|\mathbb{B}_{23}(j;5,3)|=\displaystyle\sum_{u=4}^{j-6}\displaystyle\sum_{r=0}^{u-4}\sum_{l=1}^{j-u-5}\sum_{i=1}^{j-u-l-4}P_{u-r-2}^{2}P^{2}_{j-u-l-i-2}$.
\end{corollary}

\begin{proposition}\label{Bj53524}
For $j\geq 10$, $\displaystyle|\mathbb{B}_{24}(j;5,3)|=\displaystyle\sum_{p=7}^{j-3}\sum_{q=1}^{p-6}\sum_{l=2}^{p-q-4}(l-1)P^{2}_{p-q-l-2} P_{j-p-1}^{2}$.
\end{proposition}
\begin{proof}
Let $\textbf{B}\in \mathbb{B}_{24}(j;5,3)$. Let $0<a<b<c<1$ be the reducible elements of $\textbf{B}$. As $B_{24}$(see Figure I) is the basic block associated to $\textbf{B}$, by Theorem \ref{redb}, $Red(B_{24})=Red(\textbf{B})$ and $\eta(B_{24})=\eta(\textbf{B})=3$. Observe that an adjunct representation of $B_{24}$ is given by $B_{24}=C]_{0}^{c}\{c_1\}]_{a}^{b}\{c_2\}]_{c}^{1}\{c_3\}$, where $C:0\prec a\prec x\prec b\prec c \prec y\prec 1$ is a $7$-chain. Also by Corollary \ref{maxchain}, $\textbf{B}$ has an adjunct representation $\textbf{B}=C_0]_{0}^{c}C_1]_{a}^{b}C_2]_{c}^{1}C_3$, where $C_0$ is a maximal chain containing all the reducible element of $\textbf{B}$, and $C_1, C_2, C_3$ are chains.

Now observe that $\textbf{B}=\textbf{B}'\circ \textbf{B}''$, where $\textbf{B}'=\textbf{B}\cap[0,c]\in \mathscr{B}(p;4,2,4)$ for $p\geq 7$ and $\textbf{B}''=\textbf{B}\cap[c,1]\in \mathscr{B}(q;2,1,2)$ for $q\geq 4$ with $j=p+q-1\geq 10$. Note that $7\leq p=j-q+1\leq j-3$, since $q\geq 4$. By Corollary \ref{circ}, $|\mathbb{B}_{24}(j;5,3)|=\displaystyle\sum_{p=7}^{j-3}(|\mathscr{B}(p;4,2,4)|\times|\mathscr{B}(j-p+1;2,1,2)|)$. According to Aware and Bhavale \cite{bib17}(see Proposition 3.2), $|\mathscr{B}(p;4,2,4)|=\displaystyle\sum_{i=1}^{p-6}\sum_{l=2}^{p-i-4}(l-1)P^{2}_{p-i-l-2}$. Also according to Thakare et al. \cite{bib14}(see Theorem 3.3) $|\mathscr{B}(q;2,1,2)|=P_{q-2}^{2}$. Therefore $|\mathbb{B}_{24}(j;5,3)|=\displaystyle\sum_{p=7}^{j-3}\bigg{(}\bigg{(}\displaystyle\sum_{i=1}^{p-6}\sum_{l=2}^{p-i-4}(l-1)P^{2}_{p-i-l-2}\bigg{)}\times P_{j-p-1}^{2}\bigg{)}$. Thus $|\mathbb{B}_{24}(j;5,3)|=\displaystyle\sum_{p=7}^{j-3}\sum_{i=1}^{p-6}\sum_{l=2}^{p-i-4}(l-1)P^{2}_{p-i-l-2} P_{j-p-1}^{2}$.
\end{proof}
Note that $\textbf{B}\in\mathbb{B}_{25}(j;5,3)$ if and only if $\textbf{B*}\in\mathbb{B}_{24}(j;5,3)$. Therefore using Proposition \ref{Bj53524}, we have the following result.
\begin{corollary}\label{Bj53525}
For $j\geq 10$, $\displaystyle|\mathbb{B}_{25}(j;5,3)|=\displaystyle\sum_{p=6}^{j-3}\sum_{i=1}^{p-6}\sum_{l=2}^{p-i-4}(l-1)P^{2}_{p-i-l-2} P_{j-p-1}^{2}$.
\end{corollary}

\begin{proposition}\label{Bj53526}
For $j\geq 10$, $\displaystyle|\mathbb{B}_{26}(j;5,3)|=\displaystyle\sum_{s=1}^{j-9}\sum_{t=1}^{j-s-8}\sum_{m=0}^{j-t-s-8}\sum_{p=4}^{j-t-s-m-4}P_{p-2}^{2}P_{j-t-s-m-p-2}^{2}$.
\end{proposition}
\begin{proof}
Let $\textbf{B}\in \mathbb{B}_{26}(j;5,3)$. Let $0<a<b<c<1$ be the reducible elements of $\textbf{B}$. As $B_{26}$(see Figure I) is the basic block associated to $\textbf{B}$, by Theorem \ref{redb}, $Red(B_{26})=Red(\textbf{B})$ and $\eta(B_{26})=\eta(\textbf{B})=3$. Observe that an adjunct representation of $B_{26}$ is given by $B_{26}=C]_{0}^{a}\{c_1\}]_{b}^{c}\{c_2\}]_{0}^{1}\{c_3\}$, where $C:0\prec x\prec a\prec b\prec y \prec c\prec 1$ is a $7$-chain. Also by Corollary \ref{maxchain}, $\textbf{B}$ has an adjunct representation $\textbf{B}=C_0]_{0}^{a}C_1]_{b}^{c}C_2]_{0}^{1}C_3$, where $C_0$ is a maximal chain containing all the reducible element of $\textbf{B}$, and $C_1, C_2, C_3$ are chains.

Now observe that $\textbf{B}=(\textbf{B}'\oplus C')]_{0}^{1}C_3$, where $\textbf{B}'=\textbf{B}\cap[0,c]\in\mathscr{B}(i;4,2,5)$ for $i\geq 8$, $C'=\textbf{B}\cap(c,1]$. Let $|C_3|=s\geq 1$ and $|C'|=t\geq 1$ with $j=i+t+s\geq 10$. Suppose $\textbf{D}=(\textbf{D}'\oplus C'')]_{0}^{1}C_3'\in \mathbb{B}_{26}(j;5,3)$, where $\textbf{D}'=\textbf{D}\cap[0,c]\in\mathscr{B}(i;4,2,5)$ for $i\geq 8$, $C''=\textbf{D}\cap(c,1]$ with $|C''|=t\geq 1$, and $C_3'$ is a chain with $|C_3'|=s\geq 1$. Then $\textbf{B}\cong \textbf{D}$ if and only if $\textbf{B}'\cong \textbf{D}'$, $C'\cong C''$, and $C_3\cong C'_3$.

 Now for fixed $t$ and $s$, there are $|\mathscr{B}(j-t-s;4,2,5)|$ maximal blocks in $\mathbb{B}_{26}(j;5,3)$ up to isomorphism. Therefore for fixed $s\geq 1$, $1\leq t=j-s-i\leq j-s-8$, since $i\geq 8$, and there are $\displaystyle\sum_{t=1}^{j-s-8}|\mathscr{B}(j-t-s;4,2,5)|$ maximal blocks in $\mathbb{B}_{26}(j;5,3)$ up to isomorphism. Further, $1\leq s=j-i-t\leq j-9$, since $i\geq 8$, $t\geq 1$, and there are $\displaystyle\sum_{s=1}^{j-9}\sum_{t=1}^{j-s-8}|\mathscr{B}(j-t-s;4,2,5)|$ maximal blocks in $\mathbb{B}_{26}(j;5,3)$ up to isomorphism. According to Aware and Bhavale \cite{bib17}(see Proposition 3.3), $|\mathscr{B}(i;4,2,5)|=\displaystyle\sum_{m=0}^{i-8}\sum_{p=4}^{i-m-4}P_{p-2}^{2}P_{i-m-p-2}^{2}$. Therefore there are $\displaystyle\sum_{s=1}^{j-9}\sum_{t=1}^{j-s-8}\bigg{(}\sum_{m=0}^{j-t-s-8}\sum_{p=4}^{j-t-s-m-4}P_{p-2}^{2}P_{j-t-s-m-p-2}^{2}\bigg{)}$ maximal blocks in $\mathbb{B}_{26}(j;5,3)$ up to isomorphism. Thus there are $\displaystyle\sum_{s=1}^{j-9}\sum_{t=1}^{j-s-8}\sum_{m=0}^{j-t-s-8}\sum_{p=4}^{j-t-s-m-4}P_{p-2}^{2}P_{j-t-s-m-p-2}^{2}$ maximal blocks in $\mathbb{B}_{26}(j;5,3)$ up to isomorphism.
\end{proof}
Note that $\textbf{B}\in\mathbb{B}_{27}(j;5,3)$ if and only if $\textbf{B*}\in\mathbb{B}_{26}(j;5,3)$. Therefore using Proposition \ref{Bj53526}, we have the following result.
\begin{corollary}\label{Bj53527}
For $j\geq 10$, $|\mathbb{B}_{27}(j;5,3)|=\displaystyle\sum_{s=1}^{j-9}\sum_{t=1}^{j-s-8}\sum_{m=0}^{j-t-s-8}\sum_{p=4}^{j-t-s-m-4}P_{p-2}^{2}P_{j-t-s-m-p-2}^{2}$.
\end{corollary}

\begin{proposition}\label{Bj53528}
For $j\geq 10$, $\displaystyle|\mathbb{B}_{28}(j;5,3)|=\displaystyle\sum_{t=1}^{j-9}\sum_{k=2}^{j-t-7}\sum_{l=4}^{j-k-t-3}(k-1)P_{l-2}^{2}P_{j-k-t-l-1}^{2}$.
\end{proposition}
\begin{proof}
Let $\textbf{B}\in \mathbb{B}_{28}(j;5,3)$. Let $0<a<b<c<1$ be the reducible elements of $\textbf{B}$. As $B_{28}$(see Figure I) is the basic block associated to $\textbf{B}$, by Theorem \ref{redb}, $Red(B_{28})=Red(\textbf{B})$ and $\eta(B_{28})=\eta(\textbf{B})=3$. Observe that an adjunct representation of $B_{28}$ is given by $B_{28}=C]_{a}^{b}\{c_1\}]_{b}^{c}\{c_2\}]_{0}^{1}\{c_3\}$, where $C:0\prec a\prec x \prec b\prec y \prec c \prec 1$ is a $7$-chain. Also by Corollary \ref{maxchain}, $\textbf{B}$ has an adjunct representation $\textbf{B}=C_0]_{a}^{b}C_1]_{b}^{c}C_2]_{0}^{1}C_3$, where $C_0$ is a maximal chain containing all the reducible elements of $\textbf{B}$, and $C_1, C_2, C_3$ are chains. 

Now observe that $\textbf{B}=(C'\oplus\textbf{B}'\oplus C'')]_{0}^{1}C_3$, where $\textbf{B}'=\textbf{B}\cap[a,c]\in\mathscr{B}(i;3,2,4)$ for $i\geq 7$, $C'=\textbf{B}\cap[0,a)$, $C''=\textbf{B}\cap(c,1]$ with $|C'|=k_1\geq 1$, $|C''|=k_2\geq 1$. Let $|C_3|=t\geq 1$, and let $k=k_1+k_2\geq 2$ with $j=i+k+t\geq 10$. Suppose $\textbf{D}=(E\oplus\textbf{D}'\oplus E')]_{0}^{1}C_3'\in \mathbb{B}_{28}(j;5,3)$, where $\textbf{D}'=\textbf{D}\cap[a,c]\in\mathscr{B}(i;3,2,4)$ for $i\geq 7$, $E=\textbf{D}\cap[0,a)$, $E'=\textbf{D}\cap(c,1]$ with $|E|=k_1\geq 1$, $|E'|=k_2\geq 1$, $C_3'$ is a chain with $|C_3'|=t\geq 1$. Then $\textbf{B}\cong \textbf{D}$ if and only if $\textbf{B}'\cong \textbf{D}'$, $C'\cong E$, $C''\cong E'$, and $C_3\cong C_3'$.

 Note that $k-2$(excluding $0$ and $1$) elements can be distributed over the chains $C'$ and $C''$ in $k-2+1=k-1$ ways. Now for fixed $k$ and $t$, there are $(k-1)|\mathscr{B}(j-k-t;3,2,4)|$ maximal blocks in $\mathbb{B}_{28}(j;5,3)$ up to isomorphism. Therefore for fixed $t\geq 1$, $2\leq k=j-i-t\leq j-t-7$, since $i\geq 7$, and there are $\displaystyle\sum_{k=2}^{j-t-7}(k-1)|\mathscr{B}(j-k-t;3,2,4)|$ maximal blocks in $\mathbb{B}_{28}(j;5,3)$ up to isomorphism. Further, $1\leq t=j-k-i\leq j-9$, since $i\geq 7$, $k\geq 2$, and there are $\displaystyle\sum_{t=1}^{j-9}\sum_{k=2}^{j-t-7}(k-1)|\mathscr{B}(j-k-t;3,2,4)|$ maximal blocks in $\mathbb{B}_{28}(j;5,3)$ up to isomorphism. According to Bhavale and Aware \cite{bib2}(see Corollary 3.10, put k=1), $|\mathscr{B}(i;3,2,4)|=\displaystyle\sum_{l=4}^{i-3}P_{l-2}^{2}P_{i-l-1}^{2}$. Therefore there are $\displaystyle\sum_{t=1}^{j-9}\sum_{k=2}^{j-t-7}\bigg{(}(k-1)\sum_{l=4}^{j-k-t-3}P_{l-2}^{2}P_{j-k-t-l-1}^{2}\bigg{)}$ maximal blocks in $\mathbb{B}_{28}(j;5,3)$ up to isomorphism. Thus there are $\displaystyle\sum_{t=1}^{j-9}\sum_{k=2}^{j-t-7}\sum_{l=4}^{j-k-t-3}(k-1)P_{l-2}^{2}P_{j-k-t-l-1}^{2}$ maximal blocks in $\mathbb{B}_{28}(j;5,3)$ up to isomorphism.
\end{proof}

\begin{theorem}\label{Bj536}
For $j\geq 10$, \\ \scriptsize{$|\mathscr{B}(j;5,3,6)|=\displaystyle\sum_{u=4}^{j-6}\displaystyle\sum_{r=0}^{u-4}\sum_{l=1}^{j-u-5}\sum_{s=1}^{j-u-l-4}4P_{u-r-2}^{2}P^{2}_{j-u-l-s-2}+\displaystyle\sum_{p=7}^{j-3}\sum_{q=1}^{p-6}\sum_{l=2}^{p-q-4}2(l-1)P^{2}_{p-q-l-2} P_{j-p-1}^{2}+$\\$\displaystyle\sum_{t=1}^{j-9}\sum_{k=2}^{j-t-7}\sum_{l=4}^{j-k-t-3}(k-1)P_{l-2}^{2}P_{j-k-t-l-1}^{2}+\displaystyle\sum_{s=1}^{j-9}\sum_{t=1}^{j-s-8}\sum_{m=0}^{j-t-s-8}\sum_{p=4}^{j-t-s-m-4}2P_{p-2}^{2}P_{j-t-s-m-p-2}^{2}$}.
\end{theorem}
\begin{proof}By (3) of Remark \ref{h326}, for $j\geq 10$, $\mathscr{B}(j;5,3,6)=\displaystyle\dot\cup_{i=20}^{28}\mathbb{B}_i(j;5,3)$. Therefore $\displaystyle|\mathscr{B}(j;5,3,6)|=\displaystyle\sum_{i=20}^{28}|\mathbb{B}_i(j;5,3)|$. Hence the proof follows from Proposition \ref{Bj53520}, Corollary \ref{Bj53521}, Proposition \ref{Bj53522}, Corollary \ref{Bj53523}, Proposition \ref{Bj53524}, Corollary \ref{Bj53525}, Proposition \ref{Bj53526}, Corollary \ref{Bj53527}, and Proposition \ref{Bj53528}. 
\end{proof}
\subsection{Counting of the class $\mathbb{B}(j;5,3,7)$}
In this subsection, we count the classes $\mathbb{B}_i(j;5,3)$ for $i=29$ and $i=30$; Consequently, we count the class $\mathscr{B}(j;5,3,7)$.
\begin{proposition}\label{Bj53529}
For $j\geq 11$, $\displaystyle|\mathbb{B}_{29}(j;5,3)|=\displaystyle\sum_{u=4}^{j-7}\sum_{r=0}^{u-4}\sum_{l=4}^{u-r-3}P_{u-r-2}^{2}P_{l-2}^{2}P_{u-r-l-1}^{2}$.
\end{proposition}
\begin{proof}
Let $\textbf{B}\in \mathbb{B}_{29}(j;5,3)$. Let $0<a<b<c<1$ be the reducible elements of $\textbf{B}$. As $B_{29}$(see Figure I) is the basic block associated to $\textbf{B}$, by Theorem \ref{redb}, $Red(B_{29})=Red(\textbf{B})$ and $\eta(B_{29})=\eta(\textbf{B})=3$. Observe that an adjunct representation of $B_{29}$ is given by $B_{29}=C]_{0}^{a}\{c_1\}]_{b}^{c}\{c_2\}]_{c}^{1}\{c_3\}$, where $C:0\prec x\prec a\prec b\prec y\prec c\prec z \prec 1$ is a $8$-chain. Also by Corollary \ref{maxchain}, $\textbf{B}$ has an adjunct representation $\textbf{B}=C_0]_{0}^{a}C_1]_{b}^{c}C_2]_{c}^{1}C_3$, where $C_0$ is a maximal chain containing all the reducible elements of $\textbf{B}$, and $C_1, C_2, C_3$ are chains. 

Now observe that $\textbf{B}=\textbf{B}'\oplus C'\oplus\textbf{B}''$, where $\textbf{B}'=\textbf{B}\cap [0,a]\in \mathscr{B}(p;2,1,2)$ with $p\geq 4$, $ \textbf{B}''=\textbf{B}\cap [b,c]\in \mathscr{B}(q;3,2,4)$ with $q\geq 7$, $C'$ is a chain with $|C'|=r\ge 0$, and $j=p+q+r\geq 11$. Suppose $\textbf{D}=\textbf{D}'\oplus C''\oplus\textbf{D}''\in \mathbb{B}_{29}(j;5,3)$, where $\textbf{D}'=\textbf{D}\cap [0,a]\in \mathscr{B}(p;2,1,2)$ with $p\geq 4$, $ \textbf{D}''=\textbf{D}\cap [b,c]\in \mathscr{B}(q;3,2,4)$ with $q\geq 7$, $C''$ is a chain with $|C''|=r\ge 0$. Then $\textbf{B}\cong \textbf{D}$ if and only if $\textbf{B}'\cong \textbf{D}'$, $\textbf{B}''\cong \textbf{D}''$, and $C'\cong C''$.

Let $L=\textbf{B}'\oplus C'\in \mathscr{L}'(u;2,1,2)$, where $u=p+r\geq 4$. Note that $|\mathscr{L}'(u;2,1,2)|=\displaystyle\sum_{r=0}^{u-4}P_{u-r-2}^{2}$(see proof of Proposition \ref{Bj53520}). Further, $4\leq u=j-q\leq j-7$, since $q\geq 7$. Therefore by Lemma \ref{oplus}, there are $\displaystyle\sum_{u=4}^{j-7}(|\mathscr{L}'(u;2,1,2)|\times|\mathscr{B}(j-u;3,2,4)|)$ maximal blocks in $\mathbb{B}_{29}(j;5,3)$ up to isomorphism. According to Bhavale and Aware \cite{bib2}(see Corollary 3.10, put k=1), $|\mathscr{B}(q;3,2,4)|=\displaystyle\sum_{l=4}^{q-3}P_{l-2}^{2}P_{q-l-1}^{2}$. Therefore there are $\displaystyle\sum_{u=4}^{j-7}\bigg{(}\displaystyle\sum_{r=0}^{u-4}P_{u-r-2}^{2}\times\displaystyle\sum_{l=4}^{u-r-3}P_{l-2}^{2}P_{u-r-l-1}^{2}\bigg{)}$ maximal blocks in $\mathbb{B}_{29}(j;5,3)$ up to isomorphism. Thus there are \\$\displaystyle\sum_{u=4}^{j-7}\sum_{r=0}^{u-4}\sum_{l=4}^{u-r-3}P_{u-r-2}^{2}P_{l-2}^{2}P_{u-r-l-1}^{2}$ maximal blocks in $\mathbb{B}_{29}(j;5,3)$ up to isomorphism.
\end{proof}

Note that $\textbf{B}\in\mathbb{B}_{30}(j;5,3)$ if and only if $\textbf{B*}\in\mathbb{B}_{29}(j;5,3)$. Therefore using Proposition \ref{Bj53529}, we have the following result.
\begin{corollary}\label{Bj53530}
For $j\geq 11$, $\displaystyle|\mathbb{B}_{30}(j;5,3)|=\displaystyle\sum_{u=4}^{j-7}\sum_{r=0}^{u-4}\sum_{l=4}^{u-r-3}P_{u-r-2}^{2}P_{l-2}^{2}P_{u-r-l-1}^{2}$.
\end{corollary}

\begin{theorem}\label{Bj537}
For $j\geq 11$, $|\mathscr{B}(j;5,3,7)|=\displaystyle\sum_{u=4}^{j-7}\sum_{r=0}^{u-4}\sum_{l=4}^{u-r-3}2P_{u-r-2}^{2}P_{l-2}^{2}P_{u-r-l-1}^{2}$.
\end{theorem}
\begin{proof}By (4) of Remark \ref{h326}, for $j\geq 11$, $\mathscr{B}(j;5,3,7)=\displaystyle\dot\cup_{i=29}^{30}\mathbb{B}_i(j;5,3)$. Therefore $\displaystyle|\mathscr{B}(j;5,3,7)|=\displaystyle\sum_{i=29}^{30}|\mathbb{B}_i(j;5,3)|$. Hence the proof follows from Proposition \ref{Bj53529} and Corollary \ref{Bj53530}. 
\end{proof}

\subsection{Counting of the class $\mathscr{L}(n;5,3)$}
Now in this subsection, we count the class $\mathscr{B}(j;5,3)$, and thereby, we count the class $\mathscr{L}(n;5,3)$. We have the following result.
\begin{theorem}\label{Bj53}
For $j\geq 8$, \scriptsize{$|\mathscr{B}(j;5,3)|=\displaystyle 2\binom{j-2}{6}+\displaystyle\sum_{s=1}^{j-7}\sum_{t=1}^{j-s-6}5\binom{j-t-s-2}{4}+\displaystyle\sum_{p=6}^{j-3}2\binom{p-2}{4}P_{j-p-1}^{2}+$\\$\displaystyle\sum_{p=0}^{j-9}\sum_{q=0}^{j-p-9}\sum_{r=0}^{j-p-q-9}\sum_{s=1}^{j-p-q-r-8}\sum_{t=1}^{j-p-q-r-s-7}2P_{j-p-q-r-s-t-5}^{2}+\displaystyle\sum_{t=1}^{j-8}\sum_{k=1}^{j-t-7}\sum_{p=1}^{j-t-k-6}\sum_{l=2}^{j-t-k-p-4}6(l-1)P^{2}_{j-t-k-p-l-2}+\displaystyle\sum_{t=1}^{j-8}\sum_{h=2}^{j-t-6}\sum_{l=1}^{j-h-t-5}\sum_{p=1}^{j-h-t-l-4}2(h-1)P^{2}_{j-h-t-l-p-2}+\displaystyle\sum_{u=4}^{j-6}\sum_{r=0}^{u-4}\sum_{l=1}^{j-u-5}\sum_{s=1}^{j-u-l-4}4P_{u-r-2}^{2}P^{2}_{j-u-l-s-2}+$\\$\displaystyle\sum_{p=7}^{j-3}\sum_{q=1}^{p-6}\sum_{l=2}^{p-q-4}2(l-1)P^{2}_{p-q-l-2} P_{j-p-1}^{2}+\displaystyle\sum_{t=1}^{j-9}\sum_{k=2}^{j-t-7}\sum_{l=4}^{j-k-t-3}(k-1)P_{l-2}^{2}P_{j-k-t-l-1}^{2}+\displaystyle\sum_{u=4}^{j-7}\sum_{r=0}^{u-4}\sum_{l=4}^{u-r-3}2P_{u-r-2}^{2}P_{l-2}^{2}P_{u-r-l-1}^{2}$\\$+\displaystyle\sum_{s=1}^{j-9}\sum_{t=1}^{j-s-8}\sum_{m=0}^{j-t-s-8}\sum_{p=4}^{j-t-s-m-4}2P_{p-2}^{2}P_{j-t-s-m-p-2}^{2}$}.
\end{theorem}
\begin{proof}
By (5) of Remark \ref{h326}, For $j\geq 8$, $\mathscr{B}(j;5,3)=\displaystyle\dot\cup_{h=4}^{7}\mathscr{B}(j;5,3,h)$. Therefore $\displaystyle|\mathscr{B}(j;5,3)|=\displaystyle\sum_{h=4}^{7}|\mathscr{B}(j;5,3,h)|$. Hence the proof follows from Theorem \ref{Bj534}, Theorem \ref{Bj535}, Theorem \ref{Bj536}, and Theorem \ref{Bj537}.
\end{proof}
Using Theorem \ref{Bj53}, we have the following main result.
\begin{theorem}\label{Mainthm}
For $n\geq 8$, $|\mathscr{L}(n;5,3)|=\displaystyle\sum_{i=0}^{n-8}(i+1)|\mathscr{B}(n-i;5,3)|$.
\end{theorem}
\begin{proof}
Let $L\in\mathscr{L}(n;5,3)$ with $n\geq 8$. Then $L=C\oplus\textbf{B} \oplus C'$, where $C$ and $C'$ are the chains with $|C|+|C'|=i\geq 0$, and $\textbf{B}\in\mathscr{B}(j;5,3)$ with $j=n-i\geq 8$. For fixed $i$, there are $|\mathscr{B}(n-i;5,3)|$ maximal blocks up to isomorphism. Note that there are $i+1$ ways to arrange $i$ elements on the chains $C$ and $C'$. Thus for $0\leq i=n-j\leq n-8,\displaystyle|\mathscr{L}(n;5,3)|=\sum_{i=0}^{n-8}(i+1)|\mathscr{B}(n-i;5,3)|$. Hence the result follows from Theorem \ref{Bj53}.
\end{proof}

\end{document}